\documentclass[10pt,a4paper,leqno]{article}

\usepackage[T1]{fontenc}
\usepackage[latin1]{inputenc}
\usepackage[francais,english]{babel}
\usepackage{amsmath}
\usepackage{amssymb}
\numberwithin{equation}{section}
\def\Q{\hbox{\bf Q}}
\def\N{\hbox{\bf N}}

\def\F{\hbox{\bf F}}

\def\Z{\hbox{\bf Z}}
\def\C{\hbox{\bf C}}

\DeclareMathOperator{\tr}{Tr}

\def\ds{\displaystyle}
\def\iy{\infty}
\def\mk{\medskip}

\def\b{\beta}
\def\g{\gamma}

\def\d{\delta}
\def\D{\Delta}

\def\l{\lambda}
\def\L{\Lambda}

\def\z{\zeta}

\def\e{\varepsilon}
\def\f{\varphi}

\def\s{\sigma}

\def\cf{\mathcal{F}}

\def\ck{\mathcal{K}}

\def\cs{\mathcal{S}}

\def\ni{\noindent}
\def\n{\noindent}

\def\Llr{\Longleftrightarrow}
\def\Lr{\Longrightarrow}

\def\numero{n$^{\text{o}}$}

\def\leq{\leqslant}
\def\geq{\geqslant}
\def\cleq{\preccurlyeq}
\def\clet{\prec}

\def\cget{\succ}

\newtheorem{theorem}{Th\'eor\`me}[section]
\newtheorem{definition}[theorem]{Definition}
\newtheorem{prop}[theorem]{Proposition}
\newtheorem{lem}[theorem]{Lemme}

\newtheorem{coro}[theorem]{Corollaire}
\newenvironment{dem}{\smallskip\noindent{\bf D\'emonstration}~:%
{\nopagebreak[0]}}%
{\nopagebreak[0]\hfill$\Box$ \medskip}%

\title{Formes modulaires modulo $2$ :
  L'ordre de nilpotence des op\'erateurs de Hecke
  (version d\'evelopp\'ee) \footnote{Ce texte a \'et\'e  r\'edig\'e  en 2012 \`a l'\'epoque de la
    publication de la Note \cite{NS1} et mis sur un site web devenu inaccessible.}}

\author{Jean-Louis NICOLAS\,\footnote{Recherche partiellement financ\'ee
par le CNRS, Institut Camille Jordan, UMR 5208.}}

\begin{document}
\maketitle

\begin{small}
\selectlanguage{english}
\medskip
\ni
{\bf Abstract.}
Let $\D= \sum_{m=0}^\iy
q^{(2m+1)^2} \in \F_2[[q]]$ be the reduction mod 2 of the $\D$ series.
A modular form $f$ modulo $2$ of level 1 is a polynomial in $\D$. If $p$ is an odd prime, then the Hecke
operator $T_p$ transforms $f$ in a modular form $T_p(f)$ which is a
polynomial in $\D$ whose degree is smaller than the degree of $f$, so
that $T_p$ is nilpotent.

The order of nilpotence of $f$ is defined as the smallest integer
$g=g(f)$ such that, for every family of $g$ odd primes $p_1,p_2,\ldots,p_g$, the relation
$T_{p_1}T_{p_2}\ldots T_{p_g}(f)=0$ holds. We show how one can compute
explicitly $g(f)$; if $f$ is a polynomial of degree $d\geqslant  1$  in $\D$, one finds that
$g(f)  < \frac 32 \sqrt d$.

\medskip

\ni

\medskip

\ni
\end{small}
\selectlanguage{francais}

\section{Introduction}\label{introd}

Soit 
\begin{equation*}
\D(q)=q\prod_{n=1}^\iy (1-q^n)^{24}=\sum_{n=1}^\iy \tau(n) q^n
\end{equation*}
o\`u $\tau$ est la fonction de Ramanujan.
Soit $k$ un entier $\geqslant  0$. On \'ecrit
\begin{equation*}
\D^k(q)=\sum_{n=k}^\iy \tau_k(n)q^n.
\end{equation*}
Les congruences connues sur $\tau(n) \pmod{2}$ (cf. \cite{SW1})
montrent que   
\begin{equation*}
\D(q)\equiv \sum_{m=0}^\iy q^{(2m+1)^2} \pmod{2},
\end{equation*}
ce qui entra\^ine
\begin{equation}\label{cnk}
n\not\equiv k\pmod{8}\quad \Lr \quad \tau_k(n)\equiv 0 \pmod{2}.
\end{equation}

Une forme modulaire de niveau 1 modulo $2$ est un polyn\^ome $f(\D)$  
\`a  coefficients dans $\F_2$ (cf. par exemple \cite{Nic,Ser5});
nous l'identifierons  \`a  une s\'erie formelle en la variable $q$,  \`a 
coefficients dans $\F_2$. Nous nous int\'eresserons principalement
aux formes paraboliques (celles dont le terme constant est 0).
 \`A  partir de maintenant (sauf mention expresse du contraire), 
toutes les s\'eries
consid\'er\'ees sont  \`a  coefficients mod 2, et nous nous 
permettrons d'\'ecrire
\begin{equation}\label{Delta}
\D=\D(q)=\sum_{m=0}^\iy q^{(2m+1)^2} \ \ \in \F_2[[q]].
\end{equation} 
Notons qu'une forme modulaire de poids $w$ pair et de niveau 1 \`a
coefficients entiers est congrue modulo $2$ \`a un polyn\^ome en $\D$
de degr\'e 
\begin{equation}\label{w12}
\leqslant w/12.
\end{equation} 

Les r\'esultats principaux expos\'es ci-dessous ont \'et\'e pr\'esent\'es dans la
Note \cite{NS1} (cf. aussi \cite{NS2}). Nous en donnons dans cet
article une d\'emonstration compl\`ete.
\section{Pr\'eliminaires}\label{formod}

\subsection{Les $\mathbf{\F_2}$-espaces vectoriels 
$\cf, \cf_1,\cf_3,\cf_5,\cf_7$}\label{parFi}

Soit $\cf$ le sous-espace 
de $ \F_2[\D]$ engendr\'e par $\D,\D^3,\D^5,\ldots$.
Compte tenu de (\ref{cnk}), on a 
\begin{equation}\label{Fi}
\cf=\cf_1\oplus\cf_3\oplus\cf_5\oplus\cf_7
\end{equation}
o\`u, pour $i\in\{1,3,5,7\}$, $\cf_i$ a pour base $\{\D^i,\D^{i+8},
\D^{i+16},\ldots \}.$ 

Puisque $\D^{2k}(q)=\D^k(q^2)$, toute forme parabolique modulo $2$, 
$f=\sum_{k\in \ck} \D^k$ (où $\ck$ est un ensemble de nombres entiers $>
0$) peut s'\'ecrire comme une somme finie 
\begin{equation}\label{fs}
f=\sum_{s\geqslant  0} f_s^{2^s} \quad {\rm avec} \quad f_s\in \cf,
\end{equation}
en posant
$$f_s=\sum_{k\in \ck,\, v_2(k)=s} \D^{k2^{-s}}.$$
Toute forme modulaire $f$ modulo $2$ non parabolique s'\'ecrit
\begin{equation}\label{1+fs}
f=1+\sum_{s\geqslant  0} f_s^{2^s} \quad {\rm avec} \quad f_s\in \cf.
\end{equation}

\subsection{Op\'erateurs de Hecke}

Soit $f(q)=\sum_{n\geqslant  0} c_n q^n$ une forme modulaire modulo $2$ et
soit $p$ un nombre premier $> 2$. L'op\'erateur de Hecke 
$T_p$ transforme $f$ 
en la forme 
\begin{equation}\label{ga2}
T_p|f=\sum_{n\geqslant  0} \g_n q^n \;\; \text{ avec } 
\g(n)=
\begin{cases}
c(pn)& \text {si $p$ ne divise pas $n$}\\
c(pn)+c(n/p)& \text {si $p$ divise  $n$.}
\end{cases}
\end{equation}

\begin{small}
\ni
[Nous \'ecrirons parfois $T_p(f)$  \`a   la place de $T_p|f$. 
On trouve \'egalement dans la litt\'erature la notation $f|T_p$; nous ne nous en 
  servirons pas.]
\end{small}

\mk
  
Si  $f$  est de degr\'e  $\leqslant k$  (comme polyn\^ome en $\D$), 
alors  $f$  est la r\'eduction mod 2 d'une forme modulaire de poids $12\,k$ 
et il en est de m\^eme de  $T_p|f$; on peut \'ecrire  $T_p|\D^k$ sous la forme
\begin{equation}\label{Tpfk1}
T_p|\D^k=\sum_{j=0}^k \mu_j \D^j,\quad \text{ avec }  \mu_j\in \F_2.
\end{equation}
Supposons maintenant $k$ impair. Les formules \eqref{cnk} et \eqref{ga2}
entra\^{\i}nent que
\begin{equation}\label{muj0}
j\not\equiv pk \pmod{8}\quad \Lr \quad \mu_j=0.
\end{equation}
En particulier, on a 
\begin{equation}\label{TpFi}
T_p(\cf_i) \subset \cf_j\quad \text{  si } \quad j \equiv pi \pmod{8}.
\end{equation}

L'op\'erateur de Hecke $T_p$ commute avec les op\'erations $f \mapsto
f^{2^s}$ de sorte que, si l'on conna\^it l'action de $T_p$ sur $\cf$,
par (\ref{fs}), on la conna\^it sur toutes les formes paraboliques.

%%%%%%%%%%%%%%%%%%%%%%%%%%%%%%%%%%%%%%%%%%%%%%%%%%%%%%%%%%%%%%%
\subsection{Nilpotence des op\'erateurs de Hecke modulo 2}\label{parnil}
%%%%%%%%%%%%%%%%%%%%%%%%%%%%%%%%%%%%%%%%%%%%%%%%%%%%%%%%%%%%%%%%%%%%%
 
Une des propri\'et\'es essentielles de l'op\'erateur de Hecke $T_p$ modulo 2 
est qu'il est nilpotent; autrement dit, dans \eqref{Tpfk1}, le coefficient
$\mu_k$ est nul. Lorsque $p\equiv 3,5,7 \pmod{8}$, cela r\'esulte de
\eqref{muj0}. Le cas $p\equiv 1 \pmod{8}$ est plus d\'elicat
(cf. \cite{Ha,OnoL,Ser5}). Nous en donnerons une autre preuve
au \S \ref{parpreuvenil}.

De la nilpotence de $T_p$, de \eqref{Tpfk1} et de \eqref{muj0}, on d\'eduit
pour tout $p$ premier  $\geqslant  3$, et tout $k$ impair positif,
\begin{equation}\label{Tpfk}
T_p|\D^k=\sum_{\substack{j\equiv pk \hspace{-2mm}\pmod{8}\\1\leqslant j \leqslant k-2}} 
\mu_j \D^j ,\quad \text{ avec } \mu_j\in \F_2.
\end{equation}

%%%%%%%%%%%%%%%%%%%%%%%%%%%%%%%%%%%%%%%%%%%%%%%%%%%%%%%%%%%%%%%
\subsection{D\'etermination de ${T_p|\D,T_p|\D^3,T_p|\D^5 \text{ et } 
T_p|\D^7}$}
%%%%%%%%%%%%%%%%%%%%%%%%%%%%%%%%%%%%%%%%%%%%%%%%%%%%%%%%%%%%%

\begin{prop}\label{propTp35}
{\rm (i)} Pour tout nombre premier $p$, on a $T_p|\D=0$.

{\rm (ii)} Si $p\equiv 3 \pmod{8}$, on a $T_p|\D^3=\D;$ sinon, 
on a $T_p|\D^3=0$. 

{\rm (iii)} Si $p\equiv 5 \pmod{8}$, on a $T_p|\D^5=\D;$ sinon, 
on a $T_p|\D^5=0$. 

{\rm (iv)} On a: 

$
\;\;    T_p|\D^7= 
\left\{
\begin{array}{llll}
   0  \quad {\rm si} \quad   p \equiv 1 \pmod{8}  & \! {\rm ou \ si } \quad  p \equiv -1 \pmod{16}\\ 
    \D^5 \  {\rm si}  \quad  p \equiv 3 \pmod{8}\\
   \D^3 \ {\rm si} \quad p \equiv 5 \pmod{8}\\
   \D \ \ {\rm si} \quad  p \equiv 7 \pmod{16}.
 \end{array}
\right. $
\end{prop}

\begin{small}
\ni
[Les points (ii) et (iii) figurent en exercice dans 
\cite[\S 6.7]{Ser1}. Le point (iv) \'etait connu de J. Oesterl\'e.] 
\end{small}

\mk
\ni
{\bf D\'emonstration}:

(i) Cela se d\'emontre par un calcul direct  \`a   partir de \eqref{Delta}
et \eqref{ga2}, ou bien en utilisant \eqref{Tpfk}.

(ii) Par \eqref{Tpfk}, on a $T_p|\D^3=0$ si $p\equiv 1,5$ ou $7 \pmod{8}$ 
tandis que, si $p\equiv 3\pmod{8}$, on a $T_p|\D^3=\mu \D$, avec $\mu=0$ 
ou $\mu=1$. \'Ecrivons $T_p|\D^3=\sum_{n\geqslant  1} \g_nq^n$. Par \eqref{ga2}, 
$\g_1=\mu$ est le coefficient de $q^p$ dans $\D^3=\D \D^2$. Par \eqref{Delta}, 
$\mu$ est donc congru mod 2 au nombre de solutions de l'\'equation diophantienne 
$x^2+2y^2=p$ avec $x$ et $y$ positifs et impairs. Comme $p\equiv 3\pmod{8}$, 
ce nombre de solutions est \'egal  \`a   $1$, ce qui donne $\mu=1$.

(iii) La d\'emonstration est analogue \`a celle de (ii). On a $T_p|\D^5=0$ 
si $p\equiv 1$ ou $3 \pmod{8}$, $T_p|\D^5=\l\D$ avec $\l=0$ ou $1$ 
si $p\equiv 5\pmod{8}$ et $T_p|\D^5=\mu\D^3$ avec $\mu=0$ ou $1$ si
$p\equiv 7\pmod{8}$. On montre que $\l=1$ en consid\'erant dans 
$T_p|\D^5=T_p|(\D\D^4)$ le coefficient de $q$; on montre que $\mu=0$ 
en consid\'erant celui de $q^3$.

(iv) D\'emontrons d'abord que l'on a 
\begin{equation}\label{congf7}
\D^7(q)\equiv \sum_{n\equiv 7 \; (\bmod 8)}
\frac{\s_1(n)}{8}\; q^n\pmod{2}
\end{equation}
o\`u $\s_1(n)$ d\'esigne la somme des diviseurs de $n$.

Le nombre de fa\c cons d'\'ecrire $n\equiv 7\pmod{8}$ comme somme 
$x^2+y^2+z^2+t^2$ avec
$x,y,z,t\in\Z$ est \'egal \`a
$8\s_1(n)$ (cf. \cite[Th. 386]{HW}) ; dans une telle 
repr\'esentation, aucun des quatre nombres $x,y,z,t$ n'est nul, car $n$ 
n'est pas somme de trois carr\'es. Le nombre de solutions positives de
$x^2+y^2+z^2+t^2=n$ est donc $\frac{\s_1(n)}{2}$. Maintenant, l'un 
des quatre nombres $x,y,z,t$ est congru \`a $2\pmod{4}$. Il y a donc 
$\frac{\s_1(n)}{8}$ solutions avec $x$ pair et $y,z,t$ impair; ce qui 
par \eqref{Delta}, et par $\D^7=\D^4 \D \D \D$,
d\'emontre \eqref{congf7}.

Soit $p\equiv 1\pmod{8}$. Par \eqref{Tpfk}, on a $T_p(\D^7)=0$. 

Soit $p\equiv 3\pmod{8}$. Par \eqref{Tpfk}, on a $T_p(\D^7)=0$ ou
$\D^5$. Le coefficient de $q^5$ dans $T_p(\D^7)$ est, par la formule
\eqref{congf7}, congru modulo 2 \`a  
$$\frac{\s_1(5p)}{8}=\frac{6(p+1)}{8}\equiv 1 \pmod{2}.$$

Soit $p\equiv 5\pmod{8}$. On a $T_p(\D^7)=0$ ou $\D^3$. Le coefficient de 
$q^3$ dans $T_p(\D^7)$ est, par la formule \eqref{congf7}, congru modulo 2 \`a 
$$\frac{\s_1(3p)}{8}=\frac{4(p+1)}{8}\equiv 1 \pmod{2}.$$

Soit $p\equiv 7\pmod{8}$. On a $T_p(\D^7)=0$ ou $\D$. Le coefficient de 
$q$ dans $T_p(\D^7)$ est, par la formule \eqref{congf7}, congru modulo 2 \`a

$\ds \hspace{15mm} \frac{\s_1(p)}{8}=\frac{p+1}{8}\equiv
   \begin{cases}
      1 \pmod{2}.& \text{ si } p\equiv 7\pmod{16}\\
      0 \pmod{2}.& \text{ si } p\equiv 15\pmod{16}.\hspace{15mm}\Box
   \end{cases} 
$

%%%%%%%%%%%%%%%%%%%%%%%%%%%%%%%%%%%%%%%%%%%%%%%%%%%%%%%%%%%%%%%
\subsection{L'ordre de nilpotence}\label{parindnil}
%%%%%%%%%%%%%%%%%%%%%%%%%%%%%%%%%%%%%%%%%%%%%%%%%%%%%%%%%%%%%%%
 
Par d\'efinition, {\it l'ordre de nilpotence} 
d'une forme modulaire $f\in \F_2[\D]$ est le plus petit entier 
$g=g(f)$ tel que, pour toute suite de  $g$  nombres premiers impairs  
$p_1,p_2,\ldots,p_g$, on ait 
$T_{p_1}T_{p_2}\ldots T_{p_g}|f=0$. 

\mk

\begin {small}
\ni 
[Comme les $T_p$ commutent entre eux, l'ordre dans lequel on \'ecrit 
les $T_{p_i}$ n'a pas d'importance. Noter aussi que l'on ne suppose 
pas que les  $p_i$ soient distincts.]
\end{small}

\mk
\ni
Lorsque $f$ = 0, on convient que 
\begin{equation}\label{g0}
g(f)=-\iy.
\end{equation}

Nous d\'esignerons par $g(k)=g(\D^k)$ l'ordre de nilpotence de $\D^k$.  
Comme chaque $T_p$ abaisse le degr\'e en $\D$ d'au moins 2 unit\'es, on a
\begin{equation*}
g(k) \leqslant \frac{k+1}{2}\cdot
\end{equation*}

Soit $p$ un nombre premier impair; il r\'esulte de la d\'efinition de 
l'ordre de nilpotence d'une forme modulaire $f\in\cf$ que l'on a 
\begin{equation}\label{gTpf}
g(f) \geqslant  g(T_p|f) +1.
\end{equation}

Si $f_1,f_2,\ldots,f_r \in \F_2[\D]$, on a
\begin{equation}\label{gf1f2}
g(f_1+f_2+\ldots+f_r) \leqslant \max(g(f_1),g(f_2),\ldots,g(f_r)).
\end{equation}

La proposition \ref{propTp35} permet de calculer les ordres de 
nilpotence des \'el\'ements de $\cf$ de degr\'e $\leqslant 7$; on a
\begin{eqnarray}
& &g(0)=-\infty,\;g(\D)=1,\; g(\D^3)= g(\D^3+\D)=2,\; 
\label{g1g3}\\ 
& &g(\D^5)=g(\D^5+\D)=g(\D^5+\D^3)=g(\D^5+\D^3+\D)=2\label{g5}\\
& &g(\D^7+a_5\D^5+a_3\D^3+a_1\D) =3 \quad \text{ pour } a_1,a_3,a_5
\in \F_2.\notag
\end{eqnarray}

%%%%%%%%%%%%%%%%%%%%%%%%%%%%%%%%%%%%%%%%%%%%%%%%%%%%%
%%%%%%%%%%%%%%%%%%%%%%%%%%%%%%%%%%%%%%%%%%%%%%%%%%%%%%%
\section{Calcul des $T_p|\D^k :$ une r\'ecurrence lin\'eaire}\label{parThPol}
%%%%%%%%%%%%%%%%%%%%%%%%%%%%%%%%%%%%%%%%%%%%%%%%%%%
%%%%%%%%%%%%%%%%%%%%%%%%%%%%%%%%%%%%%%%%%%%%%%%%%%

Soit $p$ un nombre premier $> 2$.
\begin{theorem}\label{thHecPol}
\hspace{-2mm}{\bf .} 
Il existe un unique polyn\^ome  sym\'etrique $F_p(X,Y) \in \F_2[X,Y]$,
\begin{equation}\label{FpYX}
F_p(X,Y)=Y^{p+1}+s_1(X) Y^p +\ldots +s_p(X) Y +s_{p+1}(X)
\end{equation}
de degr\'e $p+1$ tel que
\begin{equation}\label{recTpDk}
T_p(\D^k) = \sum_{r=1}^{p+1} s_r (\D) \ T_p (\D^{k-r})
\end{equation}
pour tout $k \geqslant  p+1$.
De plus, pour $1 \leqslant r \leqslant p+1$, $s_r(X)$ est une somme de mon\^omes en
$X$ dont les degr\'es sont congrus  \`a  $pr$ modulo $8$ et sont $\leqslant r$.

On a
\begin{equation}\label{F3XD}
F_3(X,Y)=Y^{4}+X Y +X^4
\end{equation}  
et
\begin{equation}\label{F5XD}
F_5(X,Y)=Y^{6}+X^2 Y^4 +X^4 Y^2+X Y+X^6.
\end{equation}
\end{theorem}

\ni
{\bf D\'emonstration} :
Consid\'erons d'abord, en caract\'eristique $0$, une forme modulaire
complexe $f(q)=\sum_{n=0}^\iy a_n q^n$ de poids $k$ et de niveau $1$ et $\zeta$
une racine primitive  $p$-i\`eme de l'unit\'e, auxquelles on associe
$p+1$ s\'eries de la mani\`ere suivante:

$\ds \bullet\;\; f_0=f_0(q)=p^k f(q^p)=p^k \sum_{n=0} ^\iy a_n q^{pn},$ 

\ni
et, pour $1\leqslant i \leqslant p$,

$\ds \bullet \;\;f_i=f_i(q)= f(\z ^i q^{1/p})=\sum_{n=0} ^\iy a_n 
\zeta^{ni} q^{n/p}$ 

\ni
vue comme s\'erie en $q^{1/p}$. Pour $r\geqslant  0$, on d\'efinit la somme de
Newton
$$N_r=\sum_{i=0}^p f_i^r.$$
\'Ecrivons $f^r=(\sum_{n=0}^\iy a_n q^n)^r=\sum_{n=0}^\iy b_n q^n$. Il
vient
\begin{eqnarray*}
f_0^r &=& p^{kr}f^r(q^p)=p^{kr}\sum_{n=0}^\iy b_n q^{pn}\\
f_i^r &=& f^r(\zeta^i q^{1/p})=\sum_{n=0}^\iy b_n \zeta^{ni} q^{n/p}
\end{eqnarray*}
d'o\`u,
$$N_r=p^{kr}\sum_{\substack{m=0\\p\mid m}} ^\iy b_{m/p} q^m+p\sum_{m=0}^\iy
b_{pm}q^m$$
et, par d\'efinition de l'op\'erateur de Hecke, 
\begin{equation}\label{Nr=p}
N_r=p\cdot T_p|f^r.
\end{equation} 
Ainsi, $N_r$ est une forme modulaire de poids $kr$ et de niveau
$1$. Pour $1\leqslant r \leqslant p+1$, on d\'efinit la $r$-i\`eme
fonction sym\'etrique
\begin{equation}\label{srf}
s_r=s_r(f)=\sum_{0\leqslant i_1 < i_2 < \ldots < i_r \leqslant p} f_{i_1}
f_{i_2}\ldots f_{i_r}.
\end{equation} 
On sait que les sommes $N_r$ sont reli\'ees aux fonctions sym\'etriques
$s_r$ par les formules dites de Newton:
\begin{equation}\label{Nrpetit}
N_r-s_1 N_{r-1}+ \ldots + (-1)^{r-1} s_{r-1}N_1+(-1)^r r s_r=0 \quad (1\leqslant r
\leqslant p+1)
\end{equation} 
et
 \begin{equation}\label{Nrgrand}
N_r-s_1 N_{r-1}+ \ldots -s_{p}N_{r-p}+s_{p+1}N_{r-p-1}=0 \quad (r
\geqslant  p+2).
\end{equation} 
Les formules \eqref{Nrpetit} permettent de calculer les $s_r$ en
fonction des sommes $N_r$ et montrent par r\'ecurrence que, pour $1 \leq
r \leqslant p+1$, $s_r$ est une forme modulaire de niveau $1$ et de poids
$kr$.

La formule \eqref{Nrgrand} montre que $T_p|f^r=N_r/p$ v\'erifie une
relation de r\'ecurrence lin\'eaire d'ordre $p+1$ dont le polyn\^ome
caract\'eristique est
\begin{equation}\label{FpfY}
F_p(f,Y)=\prod_{i=0}^p (Y-f_i)=Y^{p+1}+\sum_{r=1}^{p+1} (-1)^r
s_r(f) Y^r.
\end{equation} 

Supposons maintenant les coefficients $a_n$ entiers, autrement dit,
$f\in \Z[[q]]$. Pour $r\geqslant  0$, $T_p|f^r\in \Z[[q]]$, et par \eqref{Nr=p},
on a $N_r\in \Z[[q]]$. Pour $1\leqslant r \leqslant p+1$, en extrayant $s_r$ de l'\'equation
 \eqref{Nrpetit}, on voit par r\'ecurrence que
\begin{equation}\label{srZq}
r!\;s_r \in \Z[[q]].
\end{equation} 
Ensuite, en d\'eveloppant \eqref{srf}, on obtient
\begin{equation}\label{srAN}
s_r= \sum_{m=0}^\iy A_m q^{m/p} \quad \text{ avec }\quad A_m \in \Z[\zeta]
\end{equation} 
et, plus pr\'ecis\'ement,
\begin{eqnarray}\label{AN}
A_m &=& \sum_{1\leqslant i_1  < \ldots < i_{r-1} \leqslant p}
\sum_{\substack{n_0,n_1,\ldots,n_{r-1} \geq
    0\\p^2n_0+n_1+\ldots+n_{r-1}=m}}
p^k\,a_{n_0}a_{n_1}\ldots a_{n_{r-1}}
  \zeta^{i_1n_1+\ldots+i_{r-1}n_{r-1}}\notag\\
& & +\sum_{1\leqslant i_1 < \ldots < i_r \leqslant p}
\sum_{\substack{n_1,n_2,\ldots,n_{r} \geq
    0\\n_1+n_2+\ldots+n_{r}=m}}
a_{n_1}a_{n_2}\ldots a_{n_{r}}
  \zeta^{i_1n_1+\ldots+i_{r}n_{r}}.
\end{eqnarray}
En comparant \eqref{srZq} et \eqref{srAN}, on conclut
\begin{equation}\label{srqN}
s_r= \sum_{m=0}^\iy A_{pm}\; q^{m} \in \Z[[q]].
\end{equation} 
Supposons maintenant $f=\D$; on a $a_n=\tau(n)$
et, par \eqref{Delta}, 
$a_n\equiv 0\pmod{2}$ lorsque $n\not\equiv 1 \pmod{8}$. Pour calculer
$A_m\; \bmod\; 2$, on peut donc ajouter la condition
$n_i\equiv 1 \pmod{8}$ pour les indices $n_i$, $0\leqslant i\leqslant r$, figurant
dans la formule \eqref{AN}, ce qui entra\^ine
\begin{equation}\label{ANmod2}
A_m\equiv 0 \pmod{2} \quad \text{ si } \quad m\not\equiv r \pmod{8}.
\end{equation} 
Mais, sur $\C$, la forme modulaire $s_r$ est de poids $12\,r$ et donc
modulo $2$ est,
par \eqref{w12}, un polyn\^ome en $\D$ de degr\'e $\leq
r$. \eqref{srqN} et \eqref{ANmod2} montrent que les degr\'es des
mon\^omes de ce polyn\^ome sont $\equiv pr \pmod{8}$. La formule de
r\'ecurrence \eqref{recTpDk} s'obtient en r\'eduisant modulo $2$
les relations \eqref{Nrgrand} et \eqref{Nr=p}.

Comme $f_0=p^{12}\D(q^p)\equiv \D(q^p) \pmod{2}$, on a
\begin{equation}\label{FpD}
F_p(\D(q),\D(q^p))\equiv 0 \pmod{2}.
\end{equation} 

Puisque $s_r$ est un polyn\^ome en $\D$ de degr\'e $\leqslant r$, pour
conna\^itre la valeur de $s_r$ en fonction de $\D$,
Il suffit de calculer les termes du d\'eveloppement en s\'erie
de $s_r$  \`a  l'ordre $q^{r+1}$.

\mk

\ni
{\bf Exemple~:} $p=3$. Par \eqref{FpfY}, on a
\begin{eqnarray*}
F_3(\D,Y) &\equiv& (Y+q^3)\prod_{i=1}^3 (Y+\zeta^i q^{1/3} +q^3)\\
               &\equiv& (Y+q^3)((Y+q^3)^3+q) \equiv Y^4+qY+q^4
               \pmod{q^5}
\end{eqnarray*}
d'o\`u
$$F_3(\D,Y)\equiv Y^4+\D Y + \D^4 \pmod{2},$$
ce qui d\'emontre \eqref{F3XD}.

\n Vu (\ref{F3XD}), cela donne un proc\'ed\'e de calcul des $T_3|\D^k$; si $t$
est une ind\'etermin\'ee, on a: 
$$
\sum_{k=1}^{\infty} T_3(\D^k)t^k = \frac{ \D t^3}{1+\D^3t+\D^4t^4}\cdot
$$

\mk

On peut aussi calculer les sommes de Newton $N_r=T_p|\D^r \bmod 2$ 
et r\'esoudre les \'equations \eqref{Nrpetit} et \eqref{Nrgrand} 
pour d\'eterminer les $s_r$ modulo $2$.

\mk

\ni
{\bf Exemple~:} $p=5$. On a $N_1=N_2=N_3=N_4=N_6=N_8=N_9=N_{11}=0$,
$N_5=\D$, $N_7=\D^3$, $N_{10}=\D^2$ et les \'equations \eqref{Nrpetit}
et \eqref{Nrgrand} commen\c cant par $N_1,N_3,N_5,N_7,N_9,N_{11}$ donnent
successivement $s_1=N_1=0$, $s_3=N_3=0$, $s_5=N_5=0$,
$s_2=N_7/N_5=\D^2$, $s_4=s_2 N_7/N_5=\D^4$, $s_6=s_4N_7/N_5=\D^6$,
d'o\`u
$$F_5(\D,Y)\equiv Y^6+\D^2Y^4+\D^4Y^2+\D Y +\D^6 \pmod{2},$$
ce qui prouve \eqref{F5XD}.

La m\'ethode ci-dessus a \'et\'e programm\'ee avec succ\`es en SAGE par M.
Del\'eglise jusqu' \`a $p=257$. 
On trouvera la table des polyn\^omes $F_p(X,Y)$ pour $p \leqslant 257$  et
la m\'ethode de calcul sur le site \cite{siteweb}.

\n On peut calculer $T_5|\D^k$ \`a l'aide de la relation:
$$
\sum_{k=1}^{\infty} T_5(\D^k)t^k = 
\frac{ \D t^5}{1+\D^2t^2+\D^4t^4+\D^5t^5+\D^6t^6}\cdot
$$

\mk

Une autre fa\c con de proc\'eder est d'utiliser l'\'equation modulaire
(cf. \cite{Cox} ou \cite{Lang}). Soit la s\'erie d'Eisenstein
$Q=1+240\sum_{n=1}^\iy \s_3(n) q^n$, avec $\s_3(n)=\sum_{d\mid n}
d^3$. On d\'efinit l'invariant modulaire $j$ par
$$j=j(q)=\frac{Q^3}{\D}=\frac 1q +\sum_{n=0}^\iy c_nq^n=\frac 1q
+744+196884 q+21493760 q^2+\ldots$$
Comme $Q\equiv 1\pmod{2}$, on a
$$j\equiv 1/\D \pmod{2}.$$
On sait que, pour chaque  $p$, il existe un unique polynôme 
  $\Phi_p(X,Y)\in
\Z[X,Y]$, sym\'etrique et irr\'educ\-tible, unitaire de degr\'e $p+1$ 
en chaque variable, qui v\'erifie
$$\Phi_p(j(q),j(q^p))=0.$$
Supposons, comme ci-dessus, que $p \neq 2$. 
En remplaçant $j$ par $1/\D$, on a donc, par \eqref{FpD}
$$F_p(\D,Y)\equiv Y^{p+1} \D^{p+1}\Phi_p(1/\D,1/Y) \pmod{2}.$$
Le calcul de $\Phi_p$ a fait l'objet de nombreux articles~: Smith
(\cite{Smi}) a calcul\'e $\Phi_3$ (la valeur de $\Phi_3$ est
recopi\'ee dans \cite{Cox}), Berwick (\cite{Ber}) a calcul\'e $\Phi_5$,
Kaltofen et Yui (\cite{KY1,KY2}) ont calcul\'e $\Phi_7$ et $\Phi_{11}$ et Ito
(\cite{Ito}) a calcul\'e $\Phi_p$ pour $p\leqslant 53$.

\bigskip

\noindent {\bf Variante}. On peut pr\'esenter les constructions ci-dessus d'une
autre mani\`ere. Posons  $\L = \F_2[\D]$ et $\L_p = \L[Y]/(F_p(\D,Y))$.
Notons  $y$  l'image de  $Y$  dans  $\L_p$.
L'alg\`ebre $\L_p$ est un $\L$-module libre de rang $p+1$, de base
$\{1,y,...,y^p\}$, et la trace  $\tr\;:\;\L_p \to \L$ a la propri\'et\'e que
\begin{equation}\label{Trp} 
 \tr\,(y^k) = T_p| \D^k  
\end{equation}
pour tout  $k \geqslant 0$. 
 Une autre façon d'\'ecrire cette formule consiste \`a introduire la matrice
 $A$ de la multiplication par  $y$  dans le $\L$-module $\L_p$; c'est une
 matrice \`a coefficients dans  $\L$ (qui d\'epend de  $p$, bien entendu), 
et, pour tout  $k \geqslant 0$, on a :
 $\ds \tr \,(A^k) = T_p| \D^k$.\hfill $\Box$

%%%%%%%%%%%%%%%%%%%%%%%%%%%%%%%%%%%%%%%%%%%%%%%%%%%%%%%%%%%%%%%%
%%%%%%%%%%%%%%%%%%%%%%%%%%%%%%%%%%%%%%%%%%%%%%%%%%%%%%%%%%%%%%%
\section{Les op\'erateurs de Hecke  $T_3$ et $T_5$} \label{parP3P5}
%%%%%%%%%%%%%%%%%%%%%%%%%%%%%%%%%%%%%%%%%%%%%%%%%%%%%%%%%%%%%
%%%%%%%%%%%%%%%%%%%%%%%%%%%%%%%%%%%%%%%%%%%%%%%%%%%%%%%%%%%%%%%%%

%%%%%%%%%%%%%%%%%%%%%%%%%%%%%%%%%%%%%%%%%%%%%%%%%%%%%%%%%%%%%%%
\subsection{Les nombres  $\mathbf{n_3(k)}$, $\mathbf{n_5(k)}$ 
et $\mathbf{h(k)}$}\label{parn3n5h}
%%%%%%%%%%%%%%%%%%%%%%%%%%%%%%%%%%%%%%%%%%%%%%%%%%%%%%%%%%%%%

Soit $k$ un nombre entier $\geqslant  0$ et son \'ecriture dyadique $k=\sum_{i=0}^\iy
\b_i 2^i$ (avec $\b_i=0$ ou $1$). On appelle {\it support} de $k$ l'ensemble
\begin{equation}\label{Sk}
\cs(k)=\{2^i, i\geqslant  1 \text{ et } \b_i \neq 0\};
\end{equation}
on a 
$$
k=
\begin{cases}
1+\sum_{x\in \cs(k)} x & \text{ si $k$ est impair}\\
\sum_{x\in \cs(k)} x & \text{ si $k$ est pair.}
\end{cases}
$$
On note que $\cs(k)$ ne contient pas $1$.
On pose 
\begin{equation*}
n_3(k)=\sum_{i=0}^\iy \b_{2i+1} 2^i=
\sum_{\substack{i=1\\ i\;\text{ impair}\;}}^\iy \b_i 2^{\frac{i-1}{2}},\quad
n_5(k)=\sum_{i=0}^\iy \b_{2i+2} 2^i=
\sum_{\substack{i=1\\ i\;\text{ pair}\;}}^\iy \b_i 2^{\frac{i-2}{2}}.
\end{equation*}
Observons que l'on a, pour $i$ pair, $n_3(2^i)=0$ et, pour $i$ impair
$n_5(2^i)=0$. On a ainsi
\begin{equation}\label{n3n5}
n_3(k)=\sum_{i=1}^\iy \b_i n_3(2^i), \qquad n_5(k)=\sum_{i=1}^\iy \b_i
n_5(2^i).
\end{equation}
On pose ensuite
\begin{equation*}
h(k)=n_3(k)+n_5(k)=\sum_{i=1}^\iy \b_i 2^{\lfloor \frac{i-1}{2} \rfloor}
\end{equation*}
et l'on a
\begin{equation}\label{h}
h(k) =\sum_{i=1}^\iy \b_i h(2^i)  
\end{equation}
avec
\begin{equation}\label{h2i}
h(2^i)=2^{\lfloor\frac{i-1}{2}\rfloor}=
\begin{cases} n_3(2^i)=2^{\frac{i-1}{2}} & \text{si  $i$ est impair}\\
n_5(2^i)=2^{\frac{i-2}{2}} & \text{si $i$ est pair}.
\end{cases}
\end{equation}
Notons que l'on a pour $\ell \geqslant  0$
\begin{equation}\label{h2l+1}
n_3(2\ell+1)=n_3(2\ell),\quad n_5(2\ell+1)=n_5(2\ell),\quad 
h(2\ell+1)=h(2\ell).
\end{equation}
Nous appellerons $[n_3(k),n_5(k)]$ le {\it code} du nombre $k$
et nous \'ecrirons 
$$k\simeq [n_3(k),n_5(k)].$$ 
L'application $k \mapsto [n_3(k),n_5(k)]$ est une bijection de l'ensemble des nombres 
impairs (resp. pairs) $\geqslant  0$ sur $\N^2$.

Pour $k$ impair, 
les premi\`eres valeurs du code de $k$ sont donn\'ees dans la table ci-dessous.
\begin{equation}\label{tablen3+n5}
\begin{array}{|c|c|c|c|c|c|c|c|c|c|c|c|}
\hline
k&1&3&5&7&9&11&13&15&17&19&21\\
n_3,n_5& 0,0& 1,0& 0,1& 1,1& 2,0& 3,0& 2,1&3,1&0,2&1,2&0,3\\
h&0&1&1&2&2&3&3&4&2&3&3\\
\hline
\end{array}
\end{equation}
Notons que la classe de $k$ modulo $8$ donne la parit\'e de $n_3(k)$, $n_5(k)$
et $h(k)$:
\begin{equation}\label{tablen3n5h}
\begin{array}{|r|cccccccc|}
\hline
k \bmod 8 &0&1&2&3&4&5&6&7\\
\hline
n_3(k) \bmod 2&0&0&1&1&0&0&1&1\\
n_5(k) \bmod 2&0&0&0&0&1&1&1&1 \\
h(k) \bmod2&0&0&1&1&1&1&0&0\\
\hline
\end{array}.
\end{equation}
Si $k$ est impair, on a 
\begin{equation}\label{code2k1}
2k\simeq [1+2n_5(k),n_3(k)], \quad  \quad 
4k\simeq [2n_3(k),1+2n_5(k)]
\end{equation}
tandis que, si $k$ est pair, on a 
\begin{equation}\label{code2k}
2k\simeq [2n_5(k),n_3(k)], \quad \quad 
4k\simeq [2n_3(k),2n_5(k)].
\end{equation}
L'ordre de grandeur de $h(k)$ est $\sqrt k$, cf. \S \ref{parestimgk}.

% Il est facile de voir que les suites $n_3(k)/\sqrt k$ et $n_5(k)/\sqrt k$ sont 
% major\'ees. Avec un peu plus d'efforts, on peut d\'emontrer les in\'egalit\'es 
% valables pour tout $k$ impair positif:
% \begin{equation}\label{inegn3n5}
% 0\leqslant n_3(k) < \sqrt{\frac 32}\sqrt k, \quad 
% 0\leqslant n_5(k) < \frac{\sqrt 3}{2}\sqrt k
% \end{equation}
% et
% \begin{equation}\label{ineggk}
% \frac{1}{2}\sqrt k < 1+\frac{\sqrt{k-1}}{2} 
% \leqslant n_3(k)+n_5(k)+1 = h(k)+1 < \frac 32\sqrt k.
% \end{equation}
% Dans les in\'egalit\'es \eqref{inegn3n5} et \eqref{ineggk}, les coefficients de 
% $\sqrt k$ sont optimaux : on peut s'en persuader en consid\'erant les nombres
% $1+2(4^r-1)/3\simeq [2^r-1,0]$, $(4^{r+1}-1)/3\simeq [0,2^r-1]$,
% $4^r+1\simeq [0,2^{r-1}]$ et $4^r-1 \simeq [2^r-1,2^{r-1}-1]$.

%%%%%%%%%%%%%%%%%%%%%%%%%%%%%%%%%%%%%%%%%%%%%%%%%%%%%%%%%%%%%%%
\subsection{Propri\'et\'es de la fonction $h(k)$}\label{parproh}
%%%%%%%%%%%%%%%%%%%%%%%%%%%%%%%%%%%%%%%%%%%%%%%%%%%%%%%%%%%%%

\begin{lem}\label{lemh}
Soient $b$ un entier, $b\geqslant  1$. Posons 
\begin{equation}\label{H}
H(b)=\left(\sum_{i=1}^b h(2^i)\right)-2h(2^b).
\end{equation}
On a
\begin{equation}\label{lem12}
H(b)=2^{\lfloor b/2\rfloor}-2.
\end{equation}
\end{lem}

\ni
{\bf D\'emonstration} :
Posons $S(b)=\sum_{i=1}^b h(2^i)$. Lorsque $b$ est pair, $b=2a$,
par \eqref{h2i} il vient
$$S(b)=S(2a)=1+1+2+2+\ldots + 2^{a-1}+ 2^{a-1}=2^{a+1}-2$$
et
$$H(2a)=S(2a)-2h(2^{2a})=2^{a+1}-2-2\cdot 2^{a-1}=2^a-2=2^{\lfloor b/2\rfloor}-2.$$
Lorsque $b$ est impair, $b=2a-1$, il vient
$$S(b)=S(2a-1)=S(2a)-h(2^{2a})=2^{a+1}-2-2^{a-1}=3\cdot 2^{a-1}-2$$
et

$\ds H(b)=S(b)-2h(2^{b})=3\cdot 2^{a-1}-2-2\cdot
2^{a-1}=2^{a-1}-2=2^{\lfloor b/2\rfloor}-2$.
\hfill $\Box$

\mk

Soient $k$ et $\ell$ deux nombres entiers $\geqslant  0$. Soient
$\cs(k)$ et $\cs(\ell)$ leurs supports (cf. \eqref{Sk}). Si l'on a
$\cs(k)\cap\cs(\ell)=\varnothing$, cela implique
\begin{equation}\label{hklS}
h(k+\ell)=
\begin{cases}
h(k)+h(\ell) & \text{ si $k$ et $\ell$ ne sont pas tous deux
  impairs}\\
h(k)+h(\ell) +1& \text{ si } k\equiv \ell\equiv 1 \pmod{4}.
\end{cases}
\end{equation}

\begin{prop}\label{proph}
(i) Si $k$ et $\ell$ ne sont pas tous deux impairs, on a
\begin{equation}\label{hkl}
h(k+\ell)\leqslant h(k)+h(\ell).
\end{equation}
Supposons qu'il y ait \'egalit\'e dans \eqref{hkl}, alors $\cs(k) \cap \cs(\ell)$ 
est vide ou ne contient que des puissances d'exposant pair de $2$;  si
de plus 
$2^a\in \cs(k) \cap \cs(\ell)$, alors $2^{a+1} \notin \cs(k) \cup \cs(\ell)$;
et, on a
\begin{equation}\label{n5kl}
n_5(k+\ell) =  n_5(k)+n_5(\ell)-\sum_{2^a\in \cs(k) \cap \cs(\ell)} 2^{a/2}.
\end{equation}

(ii) Si $k$ et $\ell$  sont impairs, on a
\begin{equation}\label{hklbis}
h(k+\ell)\leqslant h(k)+h(\ell) +1.
\end{equation}
S'il y a \'egalit\'e dans \eqref{hklbis}, alors $k\equiv \ell\equiv 1 \pmod{4}$,
$\cs(k) \cap \cs(\ell)$ 
est vide ou ne contient que des puissances d'exposant pair de $2$,   
$2^a\in \cs(k) \cap \cs(\ell)$ entra\^ine $2^{a+1} \notin \cs(k) \cup \cs(\ell)$
et l'on a, comme en (i)
\begin{equation}\label{n5klbis}
n_5(k+\ell) =  n_5(k)+n_5(\ell)-\sum_{2^a\in \cs(k) \cap \cs(\ell)} 2^{a/2}.
\end{equation}  
\end{prop}

\begin{dem}
(i) Posons les d\'eveloppements dyadiques : $k=\sum_{i=0}^\iy \b_i 2^i$, 
$\ell=\sum_{i=0}^\iy \g_i 2^i$, $k+\ell=\sum_{i=0}^\iy \d_i 2^i$. 

Si $\cs(k) \cap \cs(\ell)=\varnothing$, par \eqref{hklS}, on a
$h(k+\ell)= h(k)+h(\ell)$ et, similairement,  $n_5(k+\ell)= n_5(k)+n_5(\ell)$ ce qui prouve
\eqref{hkl} et \eqref{n5kl}.  

Si $\cs(k) \cap \cs(\ell)\neq \varnothing$, on d\'efinit par r\'ecurrence les suites
$(a_j)_{1\leqslant j \leqslant r}$ et $(b_j)_{1\leqslant j \leqslant r}$ (avec $r \geqslant  1$) :
$a_1$ est le plus petit nombre tel que $\b_{a_1}=\g_{a_1}=1$ 
(donc $a_1 \geqslant  1$, car $k$ et $\ell$ ne sont pas tous deux impairs)
 et $b_1$ est le plus petit nombre v\'erifiant $b_1 > a_1$ et 
$\b_{b_1}=\g_{b_1}=0$. Lorsque $a_{j-1}$ et $b_{j-1}$ ont \'et\'e d\'etermin\'es, 
on choisit  pour $a_j$ le plus petit nombre tel que 
$a_j > b_{j-1}$ et $\b_{a_j}=\g_{a_j}=1$ 
et pour $b_j$ le plus petit nombre v\'erifiant $b_j > a_j$ et 
$\b_{a_j}=\g_{a_j}=0$. Le processus s'arr\^ete lorsque l'on a 
$\b_i+\g_i \leqslant 1$ pour $i > b_r$.

Notons que, si $i\notin \cup_{j=1}^r[a_j,bj]$, les chiffres $\b_i$ et
$\g_i$ s'additionnent sans retenue et l'on a
\begin{equation}\label{retenue}
i\notin \;\bigcup_{j=1}^r\;\;[a_j,bj]\quad \Lr \quad \d_i=\b_i+\g_i.
\end{equation}
Par la formule \eqref{h}, on en d\'eduit 
\begin{eqnarray*}
h(k+\ell)- h(k)-h(\ell) &=& \sum_{i=1}^\iy (\d_i-\b_i -\g_i) h(2^i)\\
&=& \sum_{j=1}^r \sum_{a_j \leqslant i \leqslant b_j} (\d_i-\b_i -\g_i) h(2^i).
\end{eqnarray*}
Pour $j$ v\'erifiant $1\leqslant j \leqslant r$, on a $\b_{a_j}=\g_{a_j}=1$,
$\d_{a_j}=0$ (donc $\d_{a_j}-\b_{a_j}-\g_{a_j}=-2$),  
$\b_{b_j}=\g_{b_j}=0$, $\d_{b_j}=1$ (donc $\d_{b_j}-\b_{b_j}-\g_{b_j}=1$),  
et, pour $a_j < i < b_j$, ou bien $\b_i+\g_i=1$ et
$\d_i=0$ ou bien $\b_i+\g_i=2$ et $\d_i=1$;
mais dans les deux cas, $\d_i-\b_i -\g_i$ est \'egal \`a  $-1$.
En utilisant la notation \eqref{H} et en appliquant \eqref{lem12}, il vient
\begin{eqnarray}\label{proph1}
 h(k+\ell)- h(k)-h(\ell) &=&
\sum_{j=1}^r \left(-2 h(2^{a_j})+h(2^{b_j})-\sum_{a_j < i < b_j} h(2^i)\right)\notag\\
\qquad =\sum_{j=1}^r \left(H(a_j)-H(b_j)\right) &=& \sum_{j=1}^r \left(2^{\lfloor a_j/2\rfloor}-
2^{\lfloor b_j/2\rfloor}\right) \leqslant 0 
\end{eqnarray} 
ce qui d\'emontre \eqref{hkl}.

\mk

Supposons  $h(k+\ell)= h(k)+h(\ell)$; il faut que, dans \eqref{proph1},
pour tout $j$, on ait $2^{\lfloor a_j/2\rfloor}=2^{\lfloor b_j/2\rfloor}$ ce qui 
impose $a_j$ pair et $b_j=a_j+1$; ainsi, les seuls \'el\'ements de 
$\cs(k) \cap \cs(\ell)$ sont $2^{a_1}, 2^{a_2}, \ldots, 2^{a_r}$. 

Par \eqref{n3n5} et \eqref{retenue}, on a
\begin{eqnarray*}
& &n_5(k+\ell)-n_5(k)-n_5(\ell) = 
\sum_{i=1}^\iy (\d_i-\b_i-\g_i) n_5(2^i)\\
&=& \sum_{j=1}^r \sum_{i=a_j}^{b_j}(\d_i-\b_i -\g_i) n_5(2^i)\\
&=& \sum_{j=1}^r (\d_{a_j}-\b_{a_j}-\g_{a_j}) n_5(2^{a_j}) +  
(\d_{a_j+1}-\b_{a_j+1}-\g_{a_j+1}) n_5(2^{a_j+1})\\
&=& \sum_{j=1}^r -2n_5(2^{a_j})= \sum_{j=1}^r -2\cdot 2^{\frac{a_j-2}{2}} =-\sum_{j=1}^r 2^{a_j/2} 
\end{eqnarray*}   
ce qui d\'emontre \eqref{n5kl} et termine la preuve de (i).

\mk
 (ii) Lorsque $k$ et $\ell$ sont impairs, par \eqref{h2l+1}, on a $h(k)=h(k-1)$, 
$h(\ell-1)=h(\ell)$, et par \eqref{hkl}, on obtient
\begin{eqnarray}
h(k+\ell)&=&h(k-1+\ell-1+2) \leqslant h(k-1+\ell-1)+h(2)\label{proph3}\\
&\leq& h(k-1)+h(\ell-1)+h(2)=h(k)+h(\ell)+1\label{proph4} 
\end{eqnarray}
ce qui prouve \eqref{hklbis}.

Supposons  $h(k+\ell)= h(k)+h(\ell) +1$; on doit avoir \'egalit\'e 
dans \eqref{proph3}; comme $\cs(2)=\{2\}$, (i)
implique  $2\notin \cs(k+\ell-2)$, ce qui entra\^ine $k+\ell-2$
multiple de $4$ et $k\equiv \ell \pmod{4}$. 
Puisque $\cs(2)\cap \cs(k+\ell-2)=\varnothing$, il vient
par \eqref{n5kl}
\begin{equation}\label{proph5}
n_5(k+\ell)= n_5(k-1+\ell-1)+n_5(2) =  n_5(k-1+\ell-1).
\end{equation}
On doit aussi avoir \'egalit\'e dans \eqref{proph4}, 
soit $h(k-1+\ell-1)=h(k-1)+h(\ell-1)$.
Comme $\cs(k)=\cs(k-1)$ et $\cs(\ell)=\cs(\ell-1)$, cela entra\^ine que
$\cs(k) \cap \cs(\ell)$ ne contient que des puissances d'exposant pair de $2$ 
(en particulier, $2\notin \cs(k) \cap \cs(\ell)$
ce qui, puisque $k\equiv \ell \pmod{4}$,
implique $k\equiv \ell\equiv 1 \pmod{4}$) et que, si 
$2^a \in \cs(k) \cap \cs(\ell)$, alors  $2^{a+1} \notin \cs(k) \cup \cs(\ell)$.
De plus, on a, par  \eqref{proph5} et \eqref{n5kl}
$$n_5(k+\ell)=n_5(k-1+\ell-1)=n_5(k-1)+n_5(\ell-1)
-\sum_{2^a \in \cs(k) \cap \cs(\ell)} 2^{a/2}$$
et, comme $n_5(k-1)=n_5(k)$ et $n_5(\ell-1)=n_5(\ell)$, cela 
prouve \eqref{n5klbis}. 
\end{dem}

\begin{coro}
Soit $k$ un nombre entier, $k\geqslant  0$. On a
\begin{equation}\label{hk+1}
h(k+1)\;
\begin{cases}
=h(k) & \text{ si $k$ est pair}\\
\leqslant h(k)+1 & \text{ si $k$ est impair}
\end{cases}
\end{equation}
\begin{equation}\label{hk+2}
h(k+2)\;
\begin{cases}
=h(k)+1 & \text{ si } k\equiv 0,1 \pmod{4}\\
\leqslant h(k) & \text{ si } k\equiv 2,3 \pmod{4}.
\end{cases}
\end{equation}
\begin{equation}\label{hk+3}
h(k+3)\leqslant h(k)+1 \qquad \text{ si } k\geqslant  0.
\end{equation}
\begin{equation}\label{hk+4}
h(k+4)\leqslant h(k)+1 \qquad \text{ si } k\geqslant  0.
\end{equation}
\end{coro}

Nous utiliserons la relation d'ordre suivante sur l'ensemble des nombres 
entiers naturels pairs (ou impairs): 
\begin{definition}\label{cord}
Si $k$ et $\ell$ ont m\^eme parit\'e, 
on dit que $\ell$ domine $k$ et on \'ecrit
\begin{equation}\label{ord}  
k\clet \ell \quad \text{ ou }\quad \ell \cget k
\end{equation}   
si l'on a
$h(k) < h(\ell)$ ou bien $h(k)=h(\ell)$ et $n_5(k)< n_5(\ell)$.
La relation 
\begin{equation}\label{ord=}  
k \cleq \ell \quad \text{ d\'efinie par $k\clet \ell$ ou $k=\ell$},
\end{equation}   
est une relation d'ordre total  
sur l'ensemble des entiers pairs (resp. impairs) $\geqslant  0$. 
\end{definition}

\begin {small}
\ni 
[Notons que 
$k\clet \ell$ n'implique pas forc\'ement $a+k \clet a+\ell$ :
exemple $a=4,k=2,\ell=4$.]
\end{small}

%%%%%%%%%%%%%%%%%%%%%%%%%%%%%%%%%%%%%%%%%%%%%%%%%%%%%%%%%%%%%%%
\subsection{La fonction $h$ pour les polyn\^omes}\label{parhPol}
%%%%%%%%%%%%%%%%%%%%%%%%%%%%%%%%%%%%%%%%%%%%%%%%%%%%%%%%%%%%%

Nous dirons qu'un polyn\^ome $P\in \F_2[x]$
est pair (resp. impair)
s'il est nul ou si tous ses exposants sont pairs (resp. s'il est non
nul et si tous ses exposants sont impairs). Dans ce paragraphe, nous
supposons que tous les polyn\^omes de $\F_2[x]$ consid\'er\'es sont soit
pairs soit impairs. \`A l'aide de la relation de domination \eqref{ord}, nous
\'ecrirons un tel polyn\^ome (non nul) sous la forme 
\begin{equation}\label{d1d2dr}
P(x)= x^{d_1}+x^{d_2}+\ldots+x^{d_r}\quad \text{ avec } \quad 
d_1 \cget d_2 \cget \ldots \cget d_r.
\end{equation}
Nous dirons que $d_1$ est l'exposant dominant de $P$ et que $x^{d_1}$ en
est le monôme dominant. 
\begin{definition}\label{defhP}
Soit $P\in \F_2[x]$, avec $P\neq 0$, \'ecrit sous la forme \eqref{d1d2dr}.
\`A  partir de $h(k)$ (cf. \eqref{h}) on d\'efinit $h(P)$ par
\begin{equation}\label{hP}
h(P)=h(d_1)=\max_{1\leqslant i \leqslant r} h(d_i).
\end{equation}
Si $P$ est le polyn\^ome nul, on pose $h(P)=-\iy$.
\end{definition}
On a
\begin{equation}\label{hP+Q}
h(P+Q)\leqslant \max(h(P),h(Q))
\end{equation}
Les formules \eqref{code2k1} et \eqref{code2k} entra\^inent
\begin{equation}\label{hP4}
h(P^4)=
\begin{cases}
2h(P)+1 & \text{ si $P$ est impair}\\
2h(P) & \text{ si $P$ est pair}.
\end{cases}
\end{equation}
De la proposition \ref{proph}, on d\'eduit:
\begin{coro}\label{corohP}
Soient deux polyn\^omes $P(X)$ et $Q(X)$  dans $\F_2[X]$; on a
\begin{equation}\label{hPQ}
h(PQ)\leqslant h(P)+h(Q)+\e
\end{equation}
o\`u $\e$ est d\'efini par
\begin{equation}\label{eps}
\e=\e(P,Q)=
\begin{cases}
0 & \text{si l'un des polyn\^omes $P$ ou $Q$ est pair}\\
1& \text{sinon.}
\end{cases}
\end{equation}
\end{coro}

\begin{small}
\ni 
[Comme le montre l'exemple des polyn\^omes $x^{4}+x^2$ et $x^{4}$, le 
mon\^ome dominant d'un produit n'est pas toujours \'egal au produit des 
mon\^omes dominants de ses facteurs.]
\end{small}

\medskip

Soient $P$ et $Q$ deux polyn\^omes non nuls \`a  coefficients 
dans $\F_2$; on suppose 
que chacun de ces polyn\^omes est soit pair soit impair et l'on d\'efinit
$\e=0$ ou $1$ par \eqref{eps}. On \'ecrit
\begin{equation}\label{lemPQPQ} 
P(x)=x^{d_1}+x^{d_2}+\ldots + x^{d_r}\quad \text{ et } \quad 
Q(x)=x^{e_1}+x^{e_2}+\ldots + x^{e_s}
\end{equation}
avec $d_1 \cget d_2 \cget \ldots \cget d_r$ et  
$e_1 \cget e_2 \cget \ldots \cget e_s$.
On d\'esigne par $w$ l'exposant dominant de $PQ$.

On d\'efinit $w'$ comme l'unique entier $\geqslant  0$ tel que  $w'\equiv d_1+e_1\pmod{2}$ et 
\begin{equation}\label{lemPQ0}
w'\simeq [n_3(d_1)+n_3(e_1)+\e, n_5(d_1)+n_5(e_1)].
\end{equation}
Soit $i$ et $j$ v\'erifiant $1\leqslant i \leqslant r$ et $1\leqslant j \leq
s$. Observons que l'on a par \eqref{hkl} ou \eqref{hklbis}
\begin{equation}\label{hdiej}
h(d_i+e_j) \leqslant h(d_i)+h(e_j)+\e \leqslant  h(d_1)+h(e_1)+\e = h(w').
\end{equation}

\begin{lem}\label{lemPQ}
(i) On a $w \cleq w'$, au sens de la relation d'ordre d\'efinie en \eqref{ord=}.

(ii) Si $\cs(d_1) \cap \cs(e_1)\neq \varnothing$, ou si $\e=1$ et l'un des
deux nombres $d_1$ ou $e_1$ est congru \`a  $3$ modulo $4$, on a
$w\clet w'$.

(iii) Si $\cs(d_1) \cap \cs(e_1) = \varnothing$, et si lorsque $\e=1$ on a 
$d_1\equiv e_1 \equiv 1 \pmod{4}$, alors le mon\^ome dominant de $PQ$ 
est le produit des mon\^omes dominants de $P$ et de $Q$ 
et l'on a $w=w'=d_1+e_1$.
\end{lem}

\begin{dem}
(i) Montrons d'abord que pour $(i,j) \neq (1,1)$, on a
\begin{equation}\label{lemPQ1}
d_i+e_j \clet w'.
\end{equation}
Si l'on a $h(d_i) < h(d_1)$ alors \eqref{hkl} ou \eqref{hklbis} entra\^ine
$$h(d_i+e_j) \leqslant h(d_i)+h(e_j)+\e < h(d_1)+h(e_1)+\e = h(w')$$
ce qui prouve \eqref{lemPQ1}; on raisonne de m\^eme si $h(e_j) < h(e_1)$.

On peut donc maintenant supposer $h(d_i)=h(d_1)$ et $h(e_j)=h(e_1)$.
Mais alors, par  \eqref{ord}, on a $n_5(d_i) < n_5(d_1)$ 
(si $i\neq 1$) ou $n_5(e_j) < n_5(e_1)$ (si $j\neq 1$).  Ensuite, 
par \eqref{hdiej}, ou bien
$h(d_i+e_j) < h(w')$ 
(ce qui entra\^ine \eqref{lemPQ1}) ou 
bien   $h(d_i+e_j) = h(w')$, ce qui, par \eqref{hdiej}, implique
$h(d_i+e_j)=h(d_i)+h(e_j)+\e$ et, par \eqref{n5kl} ou 
\eqref{n5klbis}, on a 
$$n_5(d_i+e_j) \leqslant n_5(d_i)+n_5(e_j) < n_5(d_1)+n_5(e_1) = n_5(w')$$
ce qui \`a  nouveau entra\^ine \eqref{lemPQ1}.

Pour prouver (i), il reste \`a  d\'emontrer $d_1+e_1 \cleq w'$. 
Par \eqref{hdiej}, on a $h(d_1+e_1) \leqslant h(w')$
et, ou bien $h(d_1+e_1) < h(w')$ (ce qui entra\^ine $d_1+e_1
\clet w'$) ou bien $h(d_1+e_1)=h(w')=h(d_1)+h(e_1)+\e$ et, par \eqref{n5kl} ou
\eqref{n5klbis}, on a  
$$n_5(d_1+e_1) \leqslant  n_5(d_1)+n_5(e_1) = n_5(w')$$
ce qui entra\^ine $d_1+e_1 \cleq w'$ et termine la preuve de (i).

\mk

(ii) Compte tenu de \eqref{lemPQ1}, nous devons d\'emontrer
\begin{equation}\label{lemPQ2}
d_1+e_1 \clet w'.
\end{equation}

Si $\cs(d_1) \cap \cs(e_1)$ contient une puissance d'exposant impair
de $2$ ou si $\e=1$ et que l'un des deux nombres $d_1,e_1$ est congru
\`a $3 \pmod{4}$, par la proposition \ref{proph} on a 
$h(d_1+e_1) < h(d_1)+h(e_1)+\e = h(w')$, ce qui implique
\eqref{lemPQ2}.

Si $\cs(d_1) \cap \cs(e_1)$ 
contient une  puissance de $4$, on a soit
$h(d_1+e_1) < h(d_1)+h(e_1)+\e$ (ce qui, comme pr\'ec\'edemment, implique
\eqref{lemPQ2}) soit $h(d_1+e_1) = h(d_1)+h(e_1)+\e$ et alors,
 par \eqref{n5kl} ou \eqref{n5klbis}, on a  
$$n_5(d_1+e_1) <  n_5(d_1)+n_5(e_1) = n_5(w')$$
ce qui entra\^ine \'egalement \eqref{lemPQ2} et prouve (ii).

\mk

(iii) Puisque $\cs(d_1) \cap \cs(e_1) = \varnothing$, on a
$d_1+e_1\simeq [n_3(d_1)+n_3(e_1)+\e,n_5(d_1)+n_5(e_1)]$;
on constate que le code de $d_1+e_1$ est
\'egal \`a  celui de $w'$ et, comme ces deux nombres ont m\^eme
parit\'e, ils sont \'egaux. Par \eqref{lemPQ1}, $d_1+e_1 =w'$
domine les autres exposants du produit $PQ$.
\end{dem}

%%%%%%%%%%%%%%%%%%%%%%%%%%%%%%%%%%%%%%%%%%%%%%%%%%%%%%%%%%%%%%%
\subsection{La suite des polyn\^omes $P_k^{(3)}=T_3|\D^k$} \label{pol3}
%%%%%%%%%%%%%%%%%%%%%%%%%%%%%%%%%%%%%%%%%%%%%%%%%%%%%%%%%%%%%

Utilisons les notations et les r\'esultats du \S~\!\ref{parThPol}
avec $p=3$ et posons 
$$x=\D(q),\qquad y=\D(q^3).$$ 
Appelons $K$ le corps
$\F_2(x)$ et $L$ le corps des s\'eries de Laurent en $q$, \`a savoir
$L=\F_2((q))$. Avec cette notation, on a (cf. \eqref{FpD} et \eqref{F3XD})
\begin{equation}\label{y3}
 y^4+xy+x^4=0.
\end{equation}
Le polyn\^ome $y^4+xy+x^4$ est irr\'eductible sur $K$ (cf. \cite{siteweb})
et \eqref{y3} montre que $y$ est
solution dans $L$ d'une \'equation quartique \`a coefficients dans le
corps $K$. De fa\c con plus pr\'ecise, $K(y)$ est de degr\'e $4$ sur
$K$. L'op\'erateur $T_3$ se calcule par la formule (cf. \eqref{Trp})
\begin{equation*}%\label{Tr3}
T_3(P(x))=\tr_{K(y)/K} \;P(y),
\end{equation*}
formule qui est valable pour tout polyn\^ome $P$ et m\^eme pour toute
fraction rationnelle dont le d\'enominateur ne s'annule pas en $y$, par
exemple $P(x)=x^{-k}$.

Dans ce paragraphe, pour $k\in \Z$, nous notons 
\begin{equation}\label{Trk3}
P_k(x) =P_k^{(3)} (x)=T_3(x^k)=\tr_{K(y)/K} \;y^k.
\end{equation}
On a ainsi, pour $k\in \Z$,
\begin{equation*}%\label{recP3}
P_k(x) =\tr y^k =x\tr y^{k-3}+x^4\tr y^{k-4}=x P_{k-3}(x) +x^4 P_{k-4}(x)
\end{equation*}
et
\begin{equation}\label{P3x2}
P_{2k}(x)=\tr y^{2k}=(\tr y^k)^2=P_k^2(x)=P_k(x^2).
\end{equation}
Pour $k\geqslant  0$, par \eqref{Tpfk}, il vient
\begin{equation}\label{degP3}
P_k(x)=0 \;\;\text{ ou degr\'e de }P_{k}(x)\leqslant k-2
\end{equation}
ainsi que
\begin{equation}\label{coefP3}
\text{ les exposants du polyn\^ome $P_k$ sont congrus \`a  $3k$ modulo $8$.} 
\end{equation}
On d\'eduit de \eqref{hP}, de \eqref{coefP3} et 
de \eqref{tablen3n5h} que, pour $k\geqslant  0$, 
\begin{equation}\label{hP3mod2}
h(P_k)\equiv h(3k) \equiv 
\begin{cases}
h(k)+1 \pmod{2} & \text{ si } k\not\equiv 0 \pmod{4}\\
h(k) \pmod{2} & \text{ si } k\equiv 0 \pmod{4}.
\end{cases}
\end{equation} 
Dans les deux tables ci-dessous, on trouvera les valeurs de $P_k(x)$
pour $k=0$ et pour $k$ impair, $-9 \leqslant k \leqslant 21$. Pour $k$ pair, on obtient la
valeur de $P_k(x)$ par \eqref{P3x2}.
\begin{equation}\label{tableP3}
\begin{array}{|r|cccccccccccc|}
\hline
k&0&1&3&5&7&9&11&13&15&17&19&21\\
\hline
P_k&0&0&x&0&x^5&x^3&x^9&x^7&x^{13}+x^{5}&0&x^{17}+x^{9}&x^7\\
\hline
\end{array}
\end{equation}
\begin{equation}\label{tableP3neg}
\begin{array}{|r|ccccc|}
\hline
k&1&3&5&7&9\\
\hline
P_{-k}&x^{-3}&x^{-9}&x^{-15}+x^{-7}&x^{-21}+x^{-13}&
x^{-27}+x^{-19}+x^{-11}\\
\hline
\end{array}
\end{equation}
\begin{definition}
Posons $a_0=0$, $Q_0=0$ et pour $n\geqslant  1$
\begin{equation}\label{an}
a_n=1+4+\ldots+4^{n-1}=\frac{4^n-1}{3}\simeq [0,2^{n-1}-1],\;\;
2a_n\simeq [2^n-1,0]
\end{equation}
et
\begin{equation*}
Q_n(x)=\sum_{i=1}^n x^{(a_i+3)4^{n-i}}.
\end{equation*}
\end{definition} 
Les 8 premi\`eres valeurs de $a_n$ sont 
$$\begin{array}{|r|cccccccc|}
\hline
n&0&1&2&3&4&5&6&7\\
\hline
a_n&0&1&5&21&85&341&1365&5461\\
\hline
\end{array}$$
Les supports (cf. \eqref{Sk}) de $a_n$ et de $2a_n$ sont donn\'es par
\begin{equation}\label{San}
\cs(a_n)=\cs(a_n-1)=\{4,4^2,\ldots,4^{n-1}\}, \quad
\cs(2a_n)=\{2,8,\ldots,2^{2n-1}\}.
\end{equation}
Les premiers polyn\^omes $Q_n$ sont 
\begin{eqnarray*}
Q_0\hspace{-3mm}
&=&\hspace{-3mm}0, \quad Q_1(x)=x^4, \quad
Q_2(x)=x^{16}+x^8, \quad Q_3(x)=x^{64}+x^{32}+x^{24},\\
Q_4(x)\hspace{-3mm}&=&\hspace{-3mm}
x^{256}+x^{128}+x^{96}+x^{88},\quad
Q_5(x)=x^{1024}+x^{512}+x^{384}+x^{352}+x^{344}.
\end{eqnarray*}

\begin{lem}\label{lemP30} 
Soit $n\geqslant  1$. Les polynômes $Q_n$ et $Q^2_n$ sont pairs; leurs
exposants dominants sont respectivement ${4^n}$ et ${2\cdot 4^n}$
et l'on a 
\begin{equation}\label{hQn}
h(Q_n)=h(4^n)=2^{n-1} \quad \text{ et } \quad h(Q_n^2)=h(2\cdot 4^n)=2^n.
\end{equation}
\end{lem}

\begin{dem}
Le $i$-\`eme exposant de $Q_n$, $(a_i+3)4^{n-i}$, vaut $4^n\simeq
[0,2^{n-1}]$ pour $i=1$, $2\cdot 4^{n-1}\simeq [2^{n-1},0]$ pour $i=2$
et
$$8\cdot 4^{n-i}+4^{n-i+2}+\ldots + 4^{n-2}+
4^{n-1}\simeq[2^{n-i+1},2^{n-1}-2^{n-i+1}]$$
pour $3\leqslant i \leqslant n$.

Le $i$-\`eme exposant de $Q_n^2$, $2(a_i+3)4^{n-i}$, vaut $2\cdot 4^n\simeq
[2^{n},0]$ pour $i=1$, $4^{n}\simeq [0,2^{n-1}]$ pour $i=2$
et
$$4^{n-i+2} +2\cdot 4^{n-i+2}+\ldots + 2\cdot 4^{n-2}+2\cdot 4^{n-1}
\simeq[2^n-2^{n-i+2},2^{n-i+1}]$$
pour $3\leqslant i \leqslant n$.
\end{dem}

\begin{lem}\label{lemP31}
Pour tout $n\geqslant  0$ et $k\in\Z$, on a
\begin{equation}\label{P34n+k}
P_{4^n+k}(x)=Q_n(x) P_k(x)+x^{a_n}P_{k+1}(x)
\end{equation}
et
\begin{equation}\label{P324n+k}
P_{2\cdot 4^n+k}(x)=Q_n^2(x) P_k(x)+x^{2a_n}P_{k+2}(x).
\end{equation}
\end{lem}

\begin{dem}
\`A partir de l'\'equation \eqref{y3} et des formules $a_{n+1}=4a_n+1$
et $Q_n(x)=Q_{n-1}^4(x)+x^{a_n+3}$, on d\'emontre par r\'ecurrence sur $n$
que, pour tout $n\geqslant  0$, on a
$$y^{4^n}=x^{a_n}y+Q_n(x).$$
Il vient ensuite par \eqref{Trk3}
$$P_{4^n+k}=\tr\, y^{4^n+k}=\tr\,(Q_ny^k+x^{a_n}y^{k+1})
=Q_nP_k+x^{a_n}P_{k+1},$$
ce qui d\'emontre \eqref{P34n+k}.

La preuve de \eqref{P324n+k} est similaire, \`a partir de
$y^{2\cdot4^n}=x^{2\cdot a_n}y+Q_n^2(x)$. 
\end{dem}

On d\'eduit du lemme \ref{lemP31} et des tables \eqref{tableP3}
et \eqref{tableP3neg} les valeurs particuli\`eres de $P_k$.
\begin{coro}\label{coroP34n+}
On a pour tout $n\geqslant  0$,
$$P_{4^n}=P_{4^n+1}=P_{4^n+4}=0, P_{4^n+2}=x^{a_n+1}, 
P_{4^n+3}=xQ_n, P_{4^n+5}=x^{a_n+2},$$
$$P_{2\cdot4^n}=P_{2\cdot4^n+2}=0,
P_{2\cdot4^n+1}=x^{2a_n+1}, P_{2\cdot 4^n+3}=xQ_n^2,
P_{2\cdot 4^n+4}=x^{2a_n+2},$$
$$P_{4^n-1}=x^{-3}Q_n, P_{4^n-2}=x^{-6}Q_n+x^{a_n-3}, 
P_{4^n-3}=x^{-9}Q_n+x^{a_n-6},$$
$$P_{2\cdot4^n-1}=x^{-3}Q_n^2,P_{2\cdot4^n-2}=x^{-6}Q_n^2,
P_{2\cdot4^n-3}=x^{-9}Q_n^2+x^{2a_n-3}.$$
\end{coro}

Nous avons maintenant \`a  \'etudier $h(P_k)$.

\begin{prop}\label{propP32}
Soit $k$ un nombre entier, $k\geqslant  0$, et $P_k=P_k^{(3)}$ le polyn\^ome 
introduit en \eqref{Trk3}. On a 
\begin{equation}\label{hP3k}
h(P_k) \leqslant h(k)-1
\end{equation}
o\`u la valeur de $h(P_k)$ est d\'efinie en \eqref{hP}.
\end{prop}

\ni
{\bf D\'emonstration} :
En utilisant les formules \eqref{P34n+k} et \eqref{P324n+k}, nous allons
d\'emontrer \eqref{hP3k} par r\'ecurrence sur $k$.   
\`A  l'aide des  tables \eqref{tableP3} et \eqref{tablen3+n5}
et de \eqref{P3x2} et \eqref{h2l+1},
on v\'erifie  \eqref{hP3k} pour $0 \leqslant k \leqslant 3$. Soit
maintenant $k\geqslant  4$ et supposons $h(P_j) \leqslant h(j)-1$ pour $0\leqslant j\leq
k-1$. On d\'efinit 
$n\geqslant  1$  par
\begin{equation}\label{n}
4^n\leqslant k < 4^{n+1}.
\end{equation}
Distinguons deux cas.
\begin{enumerate}
\item
{\bf Premier cas~:} $4^n \leqslant k < 2\cdot 4^n$. Introduisons $r$ par
\begin{equation*}%\label{r}
r=k-4^n.
\end{equation*}
Par \eqref{P34n+k}, on a
\begin{equation}\label{PkQnPr}
P_{k}=P_{4^n+r}=Q_n P_r+x^{a_n}P_{r+1}.
\end{equation}
On a $0\leqslant r \leqslant k-4$ et l'hypoth\`ese de r\'ecurrence s'applique
 \`a  $r$ et  \`a  $r+1$; on a ainsi 
$h(P_r)\leqslant h(r)-1$ et  $h(P_{r+1}) \leqslant h(r+1)-1$.
On a aussi $r < 2 \cdot 4^n-4^n=4^n$ et $\cs(4^n)\cap \cs(r)=\varnothing$, 
ce qui, par \eqref{hklS} et \eqref{h2i}, implique
$$h(k)=h(4^n+r) =h(4^n)+h(r)=2^{n-1}+h(r).$$
Par le lemme \ref{lemP30}, le polyn\^ome $Q_n$ est pair et,
par \eqref{hQn}, on a $h(Q_n)=2^{n-1}$. Par \eqref{hPQ}, avec
$\e(Q_n,P_r)=0$,  il vient
\begin{eqnarray}\label{QnPr}
h(Q_n P_r))  &\leq&   h(Q_n)+h(P_r) = 2^{n-1}+h(P_r)\notag\\
&\leq&  2^{n-1}+h(r)-1=h(k)-1.
\end{eqnarray}
Ensuite, par \eqref{an}, on a $h(x^{a_n})=h(a_n)=2^{n-1}-1$. 

$\bullet$ Si $k$ est impair, $r$ est aussi impair; par l'hypoth\`ese de
r\'ecurrence et par \eqref{hk+1}, on a $h(P_{r+1}) \leqslant h(r+1)-1 \leq
h(r)$. Par \eqref{coefP3}, $P_{r+1}$ est un polynôme pair, donc
$\e(x^{a_n},P_{r+1})=0$. Ainsi, par \eqref{hPQ} et \eqref{an}, il vient
\begin{equation}\label{QnPr1}
h(x^{a_n}P_{r+1})\leqslant h({a_n})+h(P_{r+1})\leqslant 2^{n-1}-1+h(r)=h(k)-1.
\end{equation}

$\bullet$ Si $k$ et $r$ sont pairs, par \eqref{an} et
\eqref{coefP3}, on a $\e(x^{a_n},P_{r+1})=1$. En compensation,
\eqref{hk+1} donne $h(r+1)=h(r)$ et \eqref{hPQ} entra\^ine 
\begin{eqnarray}\label{QnPr2}
h(x^{a_n}P_{r+1}) &\leq&  h({a_n})+h(P_{r+1}) +1\leqslant h(a_n)+h(r+1)-1+1\notag\\
& & = 2^{n-1}-1+h(r)=h(k)-1.
\end{eqnarray}

Quelle que soit la parit\'e de $k$, par \eqref{QnPr1} ou
\eqref{QnPr2}, on a donc $h(x^{a_n}P_{r+1})\leqslant h(k)-1$, ce qui, avec
\eqref{QnPr}, donne $h(P_k)\leqslant h(k)-1$, en appliquant \eqref{hP+Q}
 \`a  \eqref{PkQnPr}.

\item
{\bf Deuxi\`eme cas~:} $2\cdot 4^n \leqslant k < 4^{n+1}$. En posant
\begin{equation}\label{s}
s=k-2\cdot 4^n,
\end{equation}
par \eqref{P324n+k}, on a
\begin{equation}\label{PkQn2Ps}
P_{k}=P_{2\cdot 4^n+s}=Q_n^2 P_s+x^{2a_n}P_{s+2}.
\end{equation}
Puisque $n\geqslant  1$, on a $s\leqslant k-8$ et l'on peut appliquer
l'hypoth\`ese de r\'ecurrence  \`a  $P_s$ et $P_{s+2}$, ce qui, avec 
\eqref{hk+2}, donne
$$h(P_s)\leqslant h(s)-1 \quad \text{ et } \quad h(P_{s+2}) \leqslant h(s+2)-1
\leqslant h(s).$$
On a aussi $s < 4^{n+1}-2\cdot 4^n=2\cdot 4^n$ et 
$\cs(2\cdot 4^n) \cap \cs(s)=\varnothing$, 
ce qui, par \eqref{hklS} et \eqref{h2i}, entra\^ine
$$h(k)=h(2\cdot 4^n+s) =h(2\cdot 4^n)+h(s)=2^n+h(s).$$
$Q_n^2$ et $x^{2a_n}$ sont des polyn\^omes pairs, donc
$\e(Q_n^2,P_s)=\e(x^{2a_n},P_{s+2})=0$. Par \eqref{hPQ} et
\eqref{hQn}, il vient ensuite
\begin{equation}\label{QnPs}
h(Q_n ^2P_s)) \leqslant  h(Q_n^2)+h(P_s)=2^{n}+h(P_s)\leqslant  2^{n}+h(s)-1=
h(k)-1.
\end{equation}
Par \eqref{hPQ} et  \eqref{an}, il vient similairement
\begin{equation}\label{anPs}
h(x^{2a_n}P_{s+2})\leqslant h({2a_n})+h(P_{s+2})\leq
2^{n}-1+h(s)=h(k)-1.
\end{equation}
En appliquant \eqref{hP+Q}  \`a  \eqref{PkQn2Ps}, on obtient
$h(P_k)\leqslant h(k)-1$. \hfill $\Box$
\end{enumerate}

\begin{prop}\label{propP32i}
Soit $k$ un nombre entier, $k\geqslant 0$, et $P_k=P_k^{(3)}$ le polyn\^ome 
introduit en \eqref{Trk3}.

(i) Lorsque $k$ est multiple de $4$, on a
\begin{equation}\label{hP3k4}
h(P_k) \leqslant h(k)-2.
\end{equation}

(ii) Lorsque $k\equiv 2\pmod{4}$ et $n_5(k)=0$, on a
\begin{equation}\label{hP3k2}
h(P_k) \leqslant h(k)-3.
\end{equation}
\end{prop}

\begin{dem}
Pour prouver (i), on  observe que  \eqref{hP3k}
entra\^ine $h(P_k) \leqslant h(k)-1$; mais, 
par \eqref{hP3mod2}, on a $h(P_k)\equiv h(k) \pmod{2}$, 
ce qui implique \eqref{hP3k4}.

\mk

La preuve de (ii) se fait par r\'ecurrence sur $k$ de la m\^eme fa\c
con que la preuve de la proposition \ref{propP32}. Pour $k\leqslant 3$, il
n'y a que $k=2$ qui v\'erifie les hypoth\`eses et \eqref{hP3k2} est vraie puisque
$P_2=0$.

Supposons maintenant $k\geqslant 4$, $k\equiv 2 \pmod{4}$ et $n_5(k)=0$. On
d\'efinit $n\geqslant  1$ par \eqref{n}. La condition $n_5(k)=0$ exclut le cas
$4^n \leqslant k < 2\cdot 4^n$. On a donc $2\cdot 4^n \leqslant k < 4^{n+1}$ et
l'on d\'efinit $s$ par \eqref{s}. On a $s \leqslant  k-8$, $s=k-2\cdot
4^n\equiv 2 \pmod{4}$ et $n_5(s)=0$.
L'hypoth\`ese de r\'ecurrence appliqu\'ee  \`a 
$P_s$ donne $h(P_s)\leqslant h(s)-3$. Ensuite, comme dans \eqref{QnPs}, il vient
$$h(Q_n^2 P_s) \leqslant  h(Q_n^2)+h(P_s)=2^{n}+h(P_s)\leqslant  2^{n}+h(s)-3=
h(k)-3.$$ 

L'application de (i)  \`a  $s+2$ donne
$h(P_{s+2})\leqslant h(s+2)-2$. Et, puisque $s\equiv 2 \pmod{4}$, \eqref{hk+2}
entra\^ine $h(s+2) \leqslant h(s)$. On a ainsi 
$h(P_{s+2})\leqslant h(s)-2$. Ensuite, comme dans
\eqref{anPs}, on a
$$h(x^{2a_n}P_{s+2})\leqslant h({2a_n})+h(P_{s+2})\leq
2^{n}-1+h(s)-2=h(k)-3$$
et l'on conclut en appliquant \eqref{hP+Q}  \`a  \eqref{PkQn2Ps}.
\end{dem}

\begin{prop}\label{propP32ii}
Soit $k$ un nombre entier, $k\geqslant  0$, et $P_k=P_k^{(3)}$ le polyn\^ome 
introduit en \eqref{Trk3}.

(i) Lorsque $k$ est impair, et $n_3(k)\geqslant  1$, on a $P_k\neq 0$,
$h(P_k)=h(k)-1$ et  
l'exposant dominant (cf. \S~\!\ref{parhPol})
de $P_k$ a pour code $[n_3(k)-1, n_5(k)]$.

(ii) Lorsque $k\equiv 2 \pmod{4}$ et $n_5(k)\geqslant  1$, on a $P_k \neq
0$, $h(P_k)=h(k)-1$ et l'exposant dominant de $P_k$ a 
pour code $[n_3(k), n_5(k)-1]$.
\end{prop}

\begin{small}
\ni
[Supposons $k$ impair et $n_3(k)=0$. Si $n_5(k)=2^v$ avec $v\geqslant  0$,
$\D^k=\D^{1+2^{2v+2}}$ est une s\'erie th\^eta associ\'ee \`a $\Q(i)$
(cf. \cite[\S 9]{NS2}) et l'on a $P_k=0$. Si $n_5(k)\equiv 2^v
\pmod{2^{v+1}}$ avec $v \geqslant  0$, et $n_5(k) \neq 2^v$, on constate num\'eriquement,
pour $k< 33000$, que $P_k \neq 0$ et que l'exposant dominant de $P_k$ a pour code
$[2^{v+1}-1,n_5(k)-2^{v+1}]$.

\ni
Supposons $k$ multiple de $4$, $k=4k'$. Par \eqref{P3x2}, on a
$P_k(x)=P_{k'}(x^4)$ et, lorsque $P_{k'}\neq 0$, si l'exposant dominant de $P_{k'}$ a
pour code $[a,b]$, celui de $P_k$ a pour code $[2a,2b+1]$ si $k'$ est
impair et $[2a,2b]$ si $k'$ est pair (cf. \eqref{code2k1} et
\eqref{code2k}).

\ni
Supposons $k=2k'\equiv 2 \pmod{4}$. On a $P_k(x)=P_{k'}(x^2)$. 
Pour $k \leqslant 33000$,  on constate num\'eriquement que, lorsque
$P_k$ et $P_{k'}$ sont non nuls,
l'exposant dominant de $P_k$ est le double de l'exposant dominant de
$P_{k'}$, ce qui, via \eqref{code2k}, est prouv\'e dans la
proposition \ref{propP32ii} lorsque $n_5(k) \geqslant  1$.]
\end{small}

\mk
\ni
{\bf D\'emonstration} :
Nous allons d\'emontrer simultan\'ement (i) et (ii) par r\'ecurrence sur $k$
en suivant la m\'ethode de d\'emonstration de la proposition
\ref{propP32}.

D'abord, on v\'erifie que (i) est satisfaite pour $k=3$ et qu'il n'y a
pas d'autres valeurs de $k\leqslant 3$ satisfaisant les hypoth\`eses de la proposition.

Supposons maintenant $k\geqslant  4$ et (i) et (ii) vraies jusqu'\`a
$k-1$. On d\'efinit $n\geqslant  1$ en fonction de $k$ par \eqref{n} et l'on distingue deux cas.

\begin{enumerate}
\item
{\bf Premier cas~:} $4^n \leqslant k < 2\cdot 4^n$.
On pose $r=k-4^n \leqslant k-4$; on a par \eqref{PkQnPr}
\begin{equation}\label{PkQnPr1}
P_{k}=P_{4^n+r}=Q_n P_r+x^{a_n}P_{r+1}
\end{equation}
et l'on peut appliquer \`a $r$ et \`a $r+1$ l'hypoth\`ese de
r\'ecurrence. On a aussi $r < 2\cdot 4^n-4^n=4^n$ et $\cs(4^n)\cap
\cs(r)=\varnothing$. Par \eqref{hklS}, cela implique $h(k)=2^{n-1}+h(r)$.
On a de m\^eme $n_3(k)=n_3(r)$ et $n_5(k)=2^{n-1}+n_5(r) > 0$.

Lorsque $P_r \neq 0$, soit $m$ (resp. $M$) l'exposant
dominant de $P_r$ (resp. $Q_nP_r$). Noter que, par \eqref{coefP3}, $m$
a m\^eme parit\'e 
que $k$ et $r$ et que, par \eqref{degP3}, on a $m < r < 4^n$ ce qui
implique $\cs(4^n)\cap \cs(m)=\varnothing$.
Par le lemme \ref{lemP30}, l'exposant  dominant de $Q_n$
est ${4^n}$ et, par le lemme \ref{lemPQ} (iii), on a $M=4^n+m$. 

Lorsque $P_{r+1}\neq 0$, nous d\'esignons par ${m'}$ (resp. $w$) l'exposant
dominant de $P_{r+1}$ (resp. $x^{a_n}P_{r+1}$). Par \eqref{degP3}, on
a $m' <  r < 4^n$. Par \eqref{coefP3}, $m'$ a m\^eme parit\'e que $r+1$
et $k+1$.

On consid\`ere trois sous-cas : (a) prouve  (i) tandis que (b) et (c)
prouvent (ii).

\begin{enumerate}
\item
{\bf  $k$ et $r$ impairs et $n_3(k)\geqslant  1$.} On a
$n_3(r)=n_3(k)\geqslant  1$. Par l'hypoth\`ese de r\'ecurrence, on a 
$P_r\neq 0$, $m \simeq [n_3(r)-1,n_5(r)]$ et
\begin{equation}\label{4n+m}
M=4^n+m\simeq [n_3(r)-1,2^{n-1}+n_5(r)]= [n_3(k)-1,n_5(k)].
\end{equation}
Nous allons montrer que 
\begin{equation}\label{mondomx}
\text{ \og les exposants de $x^{a_n}P_{r+1}$ sont domin\'es par
${M}$ \fg{}}
\end{equation}
ce qui, par \eqref{PkQnPr1},  prouvera que l'exposant dominant de
$P_k$ est ${M}$ et ainsi, \'etablira (i). Si $P_{r+1}=0$,
\eqref{mondomx} est \'evidente; sinon, nous allons prouver
\begin{equation}\label{w<M}
w\clet M = 4^n+m
\end{equation}
qui entra\^ine \eqref{mondomx}.

$\bullet$ Si $r\equiv 3 \pmod{4}$ ou si $r\equiv 1 \pmod{4}$ et
$n_5(r+1)=n_5(r)=0$, par  \eqref{hP3k4}  ou \eqref{hP3k2} et par \eqref{hk+1},
il vient
$$h(P_{r+1})\leqslant h(r+1)-2 \leqslant h(r)-1$$
ce qui, par \eqref{hPQ} et \eqref{an} entra\^ine 
\begin{eqnarray*}
h(x^{a_n}P_{r+1}) &\leq& h({a_n}) + h(P_{r+1})=2^{n-1}-1+h(P_{r+1})\\ 
&\leq&  2^{n-1}-1+h(r)-1 =h(k)-2=h(4^n+m)-1
\end{eqnarray*}
ce qui prouve \eqref{w<M}.

$\bullet$ Si $r\equiv 1 \pmod{4}$ et $n_5(r+1)\geqslant  1$, l'hypoth\`ese de r\'ecurrence
donne 
$$m'\simeq [n_3(r+1),n_5(r+1)-1]= [n_3(r)+1,n_5(r)-1].$$
D\'efinissons $w'$ par 
$w'\equiv a_n+m' \equiv 1 \pmod{2}$ et
\begin{eqnarray*}
w'\hspace{-2mm} &\simeq& \hspace{-2mm}
[n_3(a_n)+n_3(m'),n_5(a_n)+n_5(m')]\\
& & =[n_3(r)+1,2^{n-1}-1+n_5(r)-1]=[n_3(k)+1,n_5(k)-2].
\end{eqnarray*}
En appliquant le lemme \ref{lemPQ} (i), on obtient $w \cleq w'$, puis,
en comparant avec \eqref{4n+m}, on voit que $w \cleq w' \clet M$, ce qui
prouve \eqref{w<M}.

\item
{\bf $k\equiv r \equiv 2 \pmod{4}$  et $n_5(r)\geqslant  1$.} On a
$n_5(r+1)=n_5(r)\geqslant  1$ et $n_3(r+1)=n_3(r)=n_3(k)\geqslant  1$,
car $2\in \cs(r)$. Par l'hypoth\`ese de
r\'ecurrence, il vient $P_r\neq 0$, $m\simeq  [n_3(r),n_5(r)-1]$,
$$4^n+m\simeq [n_3(r),2^{n-1}+n_5(r)-1]= [n_3(k),n_5(k)-1]$$
$P_{r+1}\neq 0$ (car $r+1$ est impair et $n_3(r+1)\geqslant  1$) et
$$m'\simeq [n_3(r+1)-1,n_5(r+1)]= [n_3(r)-1,n_5(r)].$$
Montrons que l'exposant dominant de $P_k$ est ${M}$
en prouvant  comme en (a)
\begin{equation}\label{w<Mb}
w\clet M = 4^n+m.
\end{equation}
Pour cela, on applique  le lemme \ref{lemPQ}  \`a $x^{a_n}P_{r+1}$. On
a  $\e=\e(x^{a_n},P_{r+1})=1$. Soit $w'$ le nombre tel que
$w'\equiv a_n+m'\equiv 0 \pmod{2}$ et 
\begin{eqnarray*}
w' &\simeq& [n_3(a_n)+n_3(m')+\e,n_5(a_n)+n_5(m')]\\
&  &=[n_3(r),2^{n-1}-1+n_5(r)]=[n_3(k),n_5(k)-1].
\end{eqnarray*}
On a donc $w'=M$. Puisque $m' < 4^n$ et $n_5(m')=n_5(r)\geqslant  1$, 
par \eqref{San}, on
a $\cs(m')\cap \cs(a_n)\neq \varnothing$ et, par le lemme \ref{lemPQ}
(ii), il en r\'esulte $w\clet w' =M$, ce qui, par \eqref{PkQnPr1},
prouve \eqref{w<Mb} et \'etablit (ii) lorsque $n_5(r) > 0$.

\item
{\bf $k\equiv r \equiv 2 \pmod{4}$  et $n_5(r)=0$.}  
On a $n_5(k)=2^{n-1} > 0$ et $n_3(r+1)=n_3(r)\geqslant  1$, 
car $2\in\cs(r)$. Par l'hypoth\`ese de
r\'ecurrence, on a $P_{r+1}\neq 0$,
$m'\simeq  [n_3(r+1)-1,0]=[n_3(r)-1,0]$ ce qui, par
\eqref{San}, entra\^ine $\cs(m')\cap \cs(a_n)=\varnothing$. On applique
alors le lemme \ref{lemPQ} (iii) \`a $x^{a_n}P_{r+1}$ avec
$\e=\e(x^{a_n},P_{r+1})=1$; par \eqref{an}, on a 
$a_n\equiv 1 \pmod{4}$ tandis que, par \eqref{coefP3}, on a 
$m'\equiv 3(r+1)\equiv 1 \pmod{4}$; et l'on a
$$w=a_n+m'\simeq [n_3(r)-1+\e,2^{n-1}-1]=
[n_3(k),n_5(k)-1].$$
Dans ce sous-cas, nous allons montrer que l'exposant dominant de
$P_k$ est ${w}$, ce qui \'etablira (ii). Pour cela, 
par \eqref{PkQnPr1}, 
il suffit de montrer que les exposants de $Q_nP_r$ sont domin\'es par 
${w}$. Mais cela r\'esulte de l'in\'egalit\'e
$$h(Q_nP_r)\leqslant h(Q_n)+h(P_r)\leqslant 2^{n-1}+h(r)-3 =h(k)-3 < h(w)$$
obtenue par \eqref{hPQ} et \eqref{hP3k2}.
\end{enumerate}

\item
{\bf Deuxi\`eme cas~: $2\cdot 4^n \leqslant k < 4^{n+1}$.}
On pose $s=k-2\cdot 4^n$ et, par \eqref{P324n+k}, on a
\begin{equation}\label{PkQn2Ps1}
P_{k}=P_{2\cdot 4^n+s}=Q_n^2 P_s+x^{2a_n}P_{s+2}.
\end{equation}
On a $s\leqslant k-8$ et on peut appliquer  l'hypoth\`ese de r\'ecurrence \`a
$s$ et \`a $s+2$. On a $s < 2\cdot 4^n$ et $\cs(2\cdot 4^n)\cap
\cs(s)= \varnothing$. Cela implique $n_3(k)=2^n+n_3(s)$, $n_5(k)=n_5(s)$ 
et $h(k)=2^n+h(s)$.

Lorsque $P_s \neq 0$, soit $m$ (resp. $M$) l'exposant
dominant de $P_s$ (resp. $Q_n^2P_s$); 
par \eqref{degP3}, on a $m < s < 2\cdot 4^n$
ce qui entra\^ine $\cs(2\cdot 4^n)\cap \cs(m)=\varnothing$.
Par le lemme \ref{lemP30}, l'exposant
dominant de $Q_n^2$ est ${2\cdot 4^n}$ et le lemme \ref{lemPQ}
(iii) donne $M=2\cdot 4^n+m$.

Nous d\'esignons par ${m'}$ (resp. $w$) l'exposant dominant de
$P_{s+2}$ (resp. $x^{2a_n}P_{s+2}$)  lorsque
$P_{s+2}$ est non nul. Par \eqref{degP3}, on a $m'\leqslant s < 2\cdot
4^n$. Notons que, par \eqref{coefP3}, $m$ et $m'$ ont m\^eme parit\'e
que $k$ et $s$.

Distinguons trois sous-cas : (a) et (b) prouvent (i) tandis que (c)
prouve (ii).
\begin{enumerate}
\item
{\bf  $k$ et $s$  impairs et $n_3(s)\geqslant  1$.} Par
l'hypoth\`ese de r\'ecurrence, on a $m\simeq [n_3(s)-1,n_5(s)]$ et
$$M=2\cdot4^n+m\simeq [2^n+n_3(s)-1,n_5(s)]= [n_3(k)-1,n_5(k)].$$
Montrons que ${M}$ est l'exposant dominant de $P_k$. Pour cela,
nous allons prouver que les exposants de $x^{2a_n}P_{s+2}$ sont
domin\'es par ${M}$.

$\bullet$ Si $s\equiv 3 \pmod{4}$, par  \eqref{hP3k} et \eqref{hk+2}, il vient
$$h(P_{s+2})\leqslant h(s+2)-1 \leqslant h(s)-1$$
tandis que, par \eqref{hPQ}, on a
$$h(x^{2a_n}P_{s+2}) \leqslant h({2a_n}) \!+ h(P_{s+2}) \leqslant \!2^n-1+h(s)-1 
\!=h(k)-2 < h(M).$$

$\bullet$ Si $s\equiv 1 \pmod{4}$, on a 
$n_3(s+2)=n_3(s)+1 \geqslant  2$ et
par l'hypoth\`ese de r\'ecurrence, on a $P_{s+2}\neq 0$,
$$m'\simeq [n_3(s+2)-1,n_5(s+2)]= [n_3(s),n_5(s)],$$
et donc $m'=s < 2\cdot 4^n$. Comme $\cs(2a_n)=\{2,8,\ldots,2^{2n-1}\}$
et que $n_3(s) > 0$, l'intersection $\cs(2a_n)\cap \cs(m')$ est non
vide. On applique le lemme \ref{lemPQ} (ii) au produit
$x^{2a_n}P_{s+2}$ en d\'efinissant $w'$ impair par $w'\simeq
[2^n-1+n_3(s),n_5(s)]$. On obtient $w \clet w' = M$ (car $w'$ et $M$  
ayant m\^eme code et m\^eme parit\'e sont \'egaux)
et les mon\^omes de $x^{2a_n}P_{s+2}$ sont
domin\'es par $x^{M}$.

\item
{\bf $k$ et $s$  impairs et $n_3(s)=0$.} On a $n_3(k)=2^n > 0$, 
$k\equiv s \equiv 1 \pmod{4}$, $n_3(s+2)=1$ et $h(s)=n_5(s)$. 
Par l'hypoth\`ese de r\'ecurrence, on a $P_{s+2}\neq 0$,
$$m'\simeq [n_3(s+2)-1,n_5(s+2)]= [0,n_5(s)]$$
et, par \eqref{San},  $\cs(2a_n)\cap \cs(m')=\varnothing$. Par le lemme
\ref{lemPQ} (iii), il vient
$$w=2a_n+m'\simeq [2^n-1,n_5(s)]=[n_3(k)-1,n_5(k)].$$
Dans ce sous-cas, l'exposant dominant de $P_k$ est ${w}$.
Pour le prouver,
il reste \`a  montrer que les exposants de
$Q_n^2P_s$ sont domin\'es par ${w}$. Si $P_s=0$, c'est clair;
sinon, soit $m$
l'exposant dominant de $P_s$. Par \eqref{hP3k}, on a
$$h(m)=n_3(m)+n_5(m)=h(P_s)\leqslant h(s)-1$$ 
ce qui entra\^ine $n_5(m) < h(s)=n_5(s)$ et l'on constate que
$$M=2\cdot 4^n+m \simeq [2^n+n_3(m),n_5(m)]\clet w=
 {2a_n+m'}\simeq [2^n-1,h(s)]$$
puisque, ou bien, $h(M) < h(w)$ 
ou bien $h(M) = h(w)$ et
$n_5(M) < n_5(w)$. 
\item
{\bf $k\equiv s \equiv 2 \pmod{4}$ et $n_5(k)=n_5(s)\geqslant  1$.}
Par l'hypoth\`ese de r\'ecurrence, on a $P_s\neq 0$, 
$m\simeq [n_3(s),n_5(s)-1]$ et
$$M=2\cdot 4^n+m \simeq [2^n+n_3(s),n_5(s)-1]=  [n_3(k),n_5(k)-1].$$
Montrons que l'exposant dominant de $P_k$ est ${M}$; 
pour cela, montrons que les exposants de $x^{2a_n}P_{s+2}$ 
sont domin\'es par ${M}$~: par 
\eqref{hP3k4}  et \eqref{hk+2}, on a $h(P_{s+2}) \leqslant h(s+2)-2
\leqslant h(s)-2$, puis, par \eqref{hPQ},
\begin{equation*}
h(x^{2a_n}P_{s+2}) \leqslant h(x^{2a_n}) + h(P_{s+2})= 2^n-1+ h(P_{s+2})
\end{equation*}

$\ds \hspace{15mm} \leqslant  2^n-1+ h(s)-2=h(k)-3 <
h(k)-1=h(M).$ \hfill $\Box$
\end{enumerate}
\end{enumerate}

%\vspace{-7mm}

%%%%%%%%%%%%%%%%%%%%%%%%%%%%%%%%%%%%%%%%%%%%%%%%%%%%%%%%%%%%%%%
\subsection{La suite des polyn\^omes $P_k^{(5)}=T_5|\D^k$}\label{pol5}
%%%%%%%%%%%%%%%%%%%%%%%%%%%%%%%%%%%%%%%%%%%%%%%%%%%%%%%%%%%%%

Dans ce paragraphe, nous \'etudions la suite de polyn\^omes $P_k^{(5)}$
en proc\'edant de fa\c con analogue, {\it mutadis mutandi}, 
\`a l'\'etude des polyn\^omes $P_k^{(3)}$ au \S\!\,\ref{pol3}. 

\mk

Posons $x=\D(q)$, $y=\D(q^5)$, $K=\F_2(x)$, $L=\F_2((q))$. Par
\eqref{FpD} et \eqref{F5XD}, il vient
\begin{equation}\label{y5}
y^6+x^2y^4+x^4y^2+xy+x^6=0.
\end{equation}
Le membre de gauche de \eqref{y5} est un polyn\^ome irr\'eductible sur
$K$ (cf. \cite{siteweb}, le corps $K(y)$ est de degr\'e $6$ sur $K$ et,
par \eqref{Trp}, on a pour $k\in \Z$,
\begin{equation*}
T_5(x^k)=\tr_{K(y)/K} \;y^k.
\end{equation*}
Dans ce paragraphe, pour $k\in \Z$, nous notons 
\begin{equation}\label{Trk5}
P_k(x)=P_k^{(5)} (x)=T_5(x^k)=\tr_{K(y)/K} \;y^k.
\end{equation}
Par \eqref{y5}, on a ainsi, 
\begin{equation*}%\label{recP5}
P_k(x)=x^2 P_{k-2}(x) +x^4 P_{k-4}(x)+ x P_{k-5}(x) +x^6 P_{k-6}(x),
\end{equation*}
\begin{equation}\label{P5x2}
P_{2k}(x)=P_k^2(x)=P_k(x^2).
\end{equation}
Pour $k\geqslant  0$, par \eqref{Tpfk}, il vient
\begin{equation}\label{coefP5}
\text{ les exposants du polyn\^ome $P_k$ sont congrus \`a  $5k$ modulo $8$.} 
\end{equation}
Lorsque $P_k\neq 0$, il r\'esulte de \eqref{coefP5} que 
le degr\'e de $P_k$ est congru \`a $k$ mod $4$; comme par \eqref{Tpfk}
ce degr\'e est $\leqslant k-2$, il vient
\begin{equation}\label{degP5}
P_k(x)=0 \;\;\text{ ou degr\'e de }P_{k}(x)\leqslant k-4.
\end{equation}
On d\'eduit de \eqref{hP}, de \eqref{coefP5} et 
de \eqref{tablen3n5h} que
\begin{equation}\label{hP5mod2}
h(P_k)\equiv h(5k)\equiv h(k)+k \pmod{2}.
\end{equation}
Dans les deux tables ci-dessous, on trouvera les valeurs de $P_k$
pour $k=0$ et pour $k$ impair, $-7\leqslant k \leqslant 21$.
Compte tenu de \eqref{P5x2}, nous ne donnons les valeurs de 
$P_k$ que pour $k$ impair. 
\begin{equation}\label{tableP5}
\begin{array}{|r|cccccccccccc|}
\hline
k=&0&1&3&5&7&9&11&13&15&17&19&21\\
\hline
P_k=&0&0&0&x&x^3&0&0&x^9&x^{11}+x^{3}&x^5&x^7&x^{17}+x^9\\
\hline
\end{array}
\end{equation}
\begin{equation}\label{tableP5neg}
\begin{array}{|r|cccc|}
\hline
k=&1&3&5&7\\
\hline
P_{-k}=&x^{-5}&x^{-15} + x^{-7}&x^{-25} + x^{-17}&x^{-35}+x^{-27}+x^{-19}\\
\hline
\end{array}
\end{equation}

\begin{definition}\label{defUVWY}
Pour $n\geqslant  0$,
les polyn\^omes $U_n,V_n,W_n$ et $Y_n$ sont respectivement d\'efinis par
$$U_0=0,\; U_1=1 \text{ et, pour }n\geqslant  2, \quad
U_n(x)=x^4U_{n-1}^4(x)+x^{12}U_{n-2}^{16}(x),$$
$$V_0=1, \;V_1=0  \text{ et, pour }n\geqslant  2, \quad
V_n(x)=x^3U_{n-1}^4(x),$$
$$W_0=0,\; W_1=0  \text{ et, pour }n\geqslant  2, \quad
W_n(x)=(x^8+x^{16})U_{n-1}^4(x)+W_{n-1}^4(x)$$
et
$$Y_0=0,\; Y_1=x^8  \text{ et, pour }n\geqslant  2, \quad
Y_n(x)=x^8U_{n}^2(x)+W_{n}^2(x).$$
\end{definition}

Pour $2\leqslant n \leqslant 5$, les valeurs des polyn\^omes  $U_n,V_n,W_n$ et $Y_n$ sont
\begin{eqnarray*}
 U_2\hspace{-3mm}
&=&\hspace{-3mm}x^4, \quad
U_3=x^{20}+x^{12}, \quad U_4=x^{84}+x^{76}+x^{52},\\
U_5\hspace{-3mm}&=&\hspace{-3mm}
x^{340}+x^{332}+x^{308}+x^{212}+x^{204},
\end{eqnarray*}
\begin{eqnarray*}
 V_2\hspace{-3mm}
&=&\hspace{-3mm}x^3, \quad
V_3=x^{19}, \quad V_4=x^{83}+x^{51},\quad
V_5=x^{339}+x^{307}+x^{211},
\end{eqnarray*}
\begin{eqnarray*}
 W_2\hspace{-3mm}&=&\hspace{-3mm} x^{16}+x^8, \quad 
W_3=x^{64}+x^{24},\quad W_4=x^{256}+x^{88}+x^{64}+x^{56}\\
W_5\hspace{-3mm}&=&\hspace{-2mm}
x^{1024}+x^{344}+x^{320}+x^{312}+x^{256}+x^{216}
\end{eqnarray*}
et
\begin{eqnarray*}
 Y_2\hspace{-3mm}
&=&\hspace{-3mm}
x^{32}, \quad Y_3=x^{128}+x^{32},\quad Y_4=x^{512}+x^{160}+x^{128}\\
Y_5\hspace{-3mm}&=&\hspace{-2mm}
x^{2048}+x^{672}+x^{640}+x^{512}+x^{416}.
\end{eqnarray*}
\begin{lem}\label{lemP50}
Soit $n\geqslant  2$. 

(i) Le polyn\^ome $U_n$ est pair et a pour degr\'e $a_n-1$
(cf. \eqref{an}). Son exposant dominant est ${a_n-1}$ et l'on a 
\begin{equation}\label{hUn}
h(xU_n)=h(U_n)=2^{n-1}-1,
\end{equation}
\begin{equation}\label{hUn2}
h(xU_n^2)=h(U_n^2)=2^n-2 \quad \text{ et } \quad h(x^3U_n^2)=2^n-1.
\end{equation}

(ii) Le polyn\^ome $V_n$ est impair et a pour degr\'e $a_n-2$.
Son exposant dominant est ${a_n-2}$ et l'on a 
\begin{equation}\label{hVn}
h(V_n)=2^{n-1}-1 \quad \text{ et } \quad h(V_n^2)=2^n-2.
\end{equation}

(iii)  Le polyn\^ome $W_n$ est pair et a pour degr\'e $4^n$.
Son exposant dominant est ${4^n}$ et l'on a 
\begin{equation}\label{hWn}
h(W_n)=2^{n-1}.
\end{equation}

(iv) Le polyn\^ome $Y_n$ est pair et a pour degr\'e $2\cdot4^n$.
Son exposant dominant est ${2^{2n+1}}$ et l'on a 
\begin{equation}\label{hYn}
h(Y_n)=2^{n}.
\end{equation}
\end{lem}

\begin{dem}
On prouve d'abord par r\'ecurrence sur $n\geqslant  2$ que le degr\'e de $U_n$ est
$a_n-1\simeq [0,2^{n-1}-1]$, ce qui implique
\begin{equation}\label{hUnmin}
h(U_n) \geqslant  2^{n-1}-1.
\end{equation}
Ensuite, on d\'emontre par r\'ecurrence sur $n\geqslant  2$ que $h(U_n) =
2^{n-1}-1$. On le v\'erifie pour $n=2$ et $n=3$. Puis, on suppose $n\geq
4$ et $h(U_\nu)=2^{\nu-1}-1$ pour $\nu < n$.

Par \eqref{hPQ} et \eqref{hP4}, il vient
$$h(x^4U_{n-1}^4) \leqslant h(x^4)+h(U_{n-1}^4)=1+2h(U_{n-1})=2^{n-1}-1$$
et, de m\^eme, $h(x^{12}U_{n-2}^{16} )\leqslant 2^{n-1}-1$. Par \eqref{hP+Q},
on en d\'eduit
$$h(U_n)\leqslant \max(h(x^4U_{n-1}^4),  h(x^{12}U_{n-2}^{16} )) \leqslant 2^{n-1}-1$$
qui, avec \eqref{hUnmin}, d\'emontre \eqref{hUn}.

Par \eqref{hPQ}, on a $h(U_n^2) \leqslant 2 h(U_n)=2^n-2$; comme le degr\'e de
$U_n^2$ est $2a_n-2$, on a $h(U_n^2) \geqslant  h(2a_n-2)=2^n-2$ et l'on
conclut que $h(U_n^2)=2^n-2$.

La preuve de $h(x^3U_n^2)=2^n-1$ est similaire, ansi que celles de
(ii), (iii) et (iv).
\end{dem}

\begin{lem}\label{lemP51}
Pour tout $n\geqslant  0$ et $k\in\Z$, on a
\begin{equation}\label{P54n+k}
P_{4^n+k}=W_n P_k+V_nP_{k+1}+U_nP_{k+4}
\end{equation}
et
\begin{equation}\label{P524n+k}
P_{2\cdot 4^n+k} = Y_nP_k+x^3U_n^2P_{k+1} +V_n^2P_{k+2} +xU_n^2P_{k+3}.
\end{equation}
\end{lem}

\begin{dem}
En multipliant \eqref{y5} par $y^2+x^2$, on obtient
$$y^8+xy^3+x^3y+x^8=0$$
et, en \'elevant au carr\'e dans $K(y)$,
$$y^{16}+x^4y^4+x^3y+x^8+x^{16}=0.$$
Puis, par r\'ecurrence sur $n\geqslant  0$, on d\'emontre
$$y^{4^n}=U_n\,y^4+V_n\, y +W_n$$
et
$$y^{2\cdot 4^n}=xU_n^2\,y^3+V_n^2\,y^2+x^3U_n^2\,y+Y_n.$$ 
Par \eqref{Trk5}, il s'ensuit que
\begin{eqnarray*}
P_{4^n+k} &=& \tr \; y^{4^n+k}=\tr \; (U_ny^{k+4}+V_n y^{k+1}+W_n
y^k)\\
&=&U_n P_{k+4} +V_n P_{k+1} + W_nP_k
\end{eqnarray*}
ce qui prouve \eqref{P54n+k}.

La preuve de \eqref{P524n+k} est similaire.
\end{dem}

\`A l'aide des tables \eqref{tableP5} et \eqref{tableP5neg}, 
on d\'eduit du lemme \ref{lemP51}:

\begin{coro}\label{coroP54n+}
Pour  $n\geqslant  0$, on a
$$P_{4^n}=P_{4^n+2}=0,\;\; P_{4^n+1}=xU_n,\;\; P_{4^n+3}=x^3U_n, 
\;\;P_{4^n+4}=xV_n,$$
$$P_{4^n-1}=x^{-5} W_n,\;\;P_{4^n-2}=x^{-10} W_n+x^{-5}V_n,$$
$$P_{2\cdot 4^n}=P_{2\cdot4^n+1}=P_{2\cdot 4^n+4}=0, \;\;P_{2\cdot 4^n+2}=x^2 U_n^2,\;\;
P_{2\cdot 4^n+3}=xV_n^2$$
et
$$P_{2\cdot 4^n-1}=x^{-5}Y_n,\;\;P_{2\cdot 4^n-2}=x^{-10}Y_n+x^{-2}U_n^2$$
\end{coro} 

Nous avons maintenant \`a  d\'eterminer $h(P_k)$.
\begin{prop}\label{propP52}
Soit $k$ un nombre entier $\geqslant  0$, $P_k=P_k^{(5)}$ le polyn\^ome  
d\'efini en \ref{Trk5}. La valeur de $h(P_k)$ est d\'efinie en \eqref{hP}.

(i) Lorsque $k$ est impair et $n_5(k)\geqslant  1$, on a 
$P_k\neq 0$, $h(P_k)=h(k)-1$ et 
l'exposant dominant (cf. \S~\!\ref{parhPol})
de $P_k$ a pour code $[n_3(k), n_5(k)-1]$. 

(ii) Lorsque $k$ est impair et $n_5(k) = 0$, on a 
\begin{equation}\label{hP5k1}
h(P_k) \leqslant h(k)-3.
\end{equation}

(iii) Lorsque $k$ est pair, on a
\begin{equation}\label{hP5k2}
h(P_k) \leqslant h(k)-2.
\end{equation}

(iv) Pour tout $k$, on a 
\begin{equation}\label{hP5k}
h(P_k) \leqslant h(k)-1.
\end{equation}
\end{prop}

\begin{small}
\ni
[Supposons $k$ impair et $n_5(k)=0$. Si $n_3(k)=2^v$ avec $v\geqslant  0$,
$\D^k=\D^{1+2^{2v+1}}$ est une s\'erie th\^eta associ\'ee \`a $\Q(\sqrt{-2})$
(cf. \cite[\S 8]{NS2}) et l'on a $P_k=0$. 
Si $n_3(k)=2^v-1$ avec $v\geqslant  0$, on a $k=(1+2^{2v+1})/3$; la forme
$\D^k$ qui, par le corollaire \ref{coroP34n+}, est \'egale \`a
$T_3|\D^{1+2^{2v+1}}$ est aussi une s\'erie th\^eta associ\'ee \`a
$\Q(\sqrt{-2})$ et l'on a $P_k=0$.
Si $n_3(k)\equiv 2^v$ ou $2^v-1\pmod{2^{v+1}}$ avec $v\geqslant  1$ et $n_3(k) >
2^v$, on constate num\'eriquement que, 
pour $k< 33000$, $P_k \neq 0$ et que l'exposant dominant de $P_k$ a
pour code $[n_3(k)-2^{v+1},2^{v}-1]$.

\ni
Supposons $k$ pair, $k=2k'$. Par \eqref{P5x2}, on a
$P_k(x)=P_{k'}(x^2)$. On constate
num\'eriquement pour $k \leqslant 33000)$ que, ou bien $P_k=P_{k'}=0$ ou bien
l'exposant dominant de $P_k$ est le double de l'exposant dominant de
$P_{k'}$.]

\end{small}

\mk
\ni
{\bf D\'emonstration} :
Le point (iv) r\'esulte de (i), (ii) et (iii). Nous allons d\'emontrer
simultan\'ement (i), (ii) et (iii) par r\'ecurrence sur $k$, en suivant le
m\^eme plan que pour la d\'emonstration des  propositions \ref{propP32}
et \ref{propP32ii}.

Lorsque $k\leqslant 15$, \`a l'aide des tables \eqref{tableP5} et
\eqref{tablen3+n5}  et de \eqref{P5x2}, on constate 
que (i), (ii) et (iii) sont v\'erifi\'ees. Rappelons que, si $P_k=0$,
$h(P_k)=-\iy$.

Soit maintenant $k\geqslant  16$ et supposons (i), (ii) et (iii) vraies
jusqu'\`a $k-1$. On d\'efinit $n\geqslant  2$ par \eqref{n} et l'on distingue
deux cas.

\begin{enumerate}
\item
{\bf Premier cas~: $4^n \leqslant k < 2\cdot 4^n$.}
On pose $r=k-4^n \leqslant k-16$; on a par \eqref{P54n+k}~:
\begin{equation}\label{P54n+r}
P_{k}=P_{4^n+r}=W_n P_r+V_nP_{r+1}+U_nP_{r+4}
\end{equation}
et l'on peut appliquer \`a $r$, $r+1$ et $r+4$ l'hypoth\`ese de
r\'ecurrence. Notons  que l'on a $r\leqslant k-4^n < 4^n$, $n_3(k)=n_3(r)$,
$n_5(k)=2^{n-1}+n_5(r) > 0$ et $h(k)=2^{n-1}+h(r)$. 

Consid\'erons trois sous-cas : (a) et (b) prouvent (i) tandis que (c)
prouve (iii); puisque $n_5(k) > 0$, (ii) ne se pr\'esente pas. 

\begin{enumerate}
\item
{\bf $k$ et $r$ impairs et $n_5(r)\geqslant  1$.} Nous allons
montrer que le mon\^ome dominant de $W_nP_r$ est aussi 
le mon\^ome dominant de $P_{k}$.

Par le lemme \ref{lemP50} (iii), 
l'exposant dominant de $W_n$ est ${4^n}$
tandis que, par l'hypoth\`ese de r\'ecurrence, $P_r$ est non nul et
l'exposant dominant $m$ de $P_r$ a pour code $[n_3(r),n_5(r)-1]$.
Par \eqref{degP5}, on a $m < r < 4^n$, ce qui entra\^ine $\cs(4^n) \cap
\cs(m)=\varnothing$;
par le lemme \ref{lemPQ} (iii), 
l'exposant dominant de $W_nP_r$ est  
$$M=4^n+m \simeq [n_3(r),2^{n-1}+n_5(r)-1] =[n_3(k),n_5(k)-1]$$
et $h(M)=h(k)-1$.

Les exposants de $V_nP_{r+1}$ sont domin\'es par $M$. 
On a en effet par 
\eqref{hPQ} (en notant que $P_{r+1}$ est pair), \eqref{hVn},
l'hypoth\`ese de r\'ecurrence  et \eqref{hk+1} 
\begin{eqnarray}\label{propP523}
h(V_nP_{r+1}) &\leq& h(V_n)+h(P_{r+1}) \leqslant 2^{n-1}-1+h(r+1)-2\notag\\
&\leq& 2^{n-1}-3+h(r)+1 = h(k)-2,
\end{eqnarray}
ce qui implique $h(V_n P_{r+1}) < h(M)$.

Lorsque $P_{r+4}\neq 0$, il reste \`a  consid\'erer le polyn\^ome
$U_nP_{r+4}$. Par le lemme \ref{lemP50} (i) et \eqref{hUn}, 
l'exposant dominant de $U_n$ est ${a_n-1}$ et $h(U_n)=h(a_n-1)=2^{n-1}-1$.

$\bullet$ {\it Supposons} $n_5(r+4)=0$; en appliquant successivement
\eqref{hPQ} (en notant que $U_n$ est pair), \eqref{hUn},
l'hypoth\`ese de r\'ecurrence et \eqref{hk+4}, on a 
\begin{eqnarray*}
h(U_nP_{r+4}) &\leq& h(U_n)+h(P_{r+4}) \leqslant 2^{n-1}-1+h(r+4)-3\\
&\leq& 2^{n-1}-4+h(r)+1 = h(k)-3 < h(M),
\end{eqnarray*}
ce qui montre que les exposants de $U_nP_{r+4}$ sont domin\'es par $M$.

$\bullet$ {\it Supposons} $n_5(r+4) > 0$; par l'hypoth\`ese de
r\'ecurrence, l'exposant dominant $m'$ 
de $P_{r+4}$ a pour code 
\begin{equation}\label{m'5}
m'\simeq [n_3(r+4), n_5(r+4)-1].
\end{equation}
D\'esignons par $w$ l'exposant dominant  de $U_nP_{r+4}$ et appliquons 
le lemme \ref{lemPQ} au produit $U_nP_{r+4}$. 
Le nombre impair $w'$ a pour code
$$w'\simeq [n_3(r+4),2^{n-1}+n_5(r+4)-2].$$

$\quad$Si $h(r+4) \leqslant h(r)$, on a $h(w') < h(M)$ donc $w' \clet M$ 
et le lemme \ref{lemPQ} (i) donne $w\cleq w' \clet M$.

$\quad$Si $h(r+4) > h(r)$, \eqref{hk+4} entra\^ine $h(r+4)=h(r)+1$ et
$h(w')=h(M)$. Distinguons deux possibilit\'es.

$\quad$Si $4\in \cs(r)$, par \eqref{n5kl}, on a
$n_5(r+4)\leqslant n_5(r)+n_5(4)-2=n_5(r)-1$, ce qui implique $w'\clet M$ 
et le lemme \ref{lemPQ} (i)  donne $w\cleq w' \clet M$.

$\quad$Si $4\notin \cs(r)$, on a $n_5(r+4)=n_5(r)+1 \geq
2$, $n_3(r+4)=n_3(r)$ et $w'=M$. Par \eqref{degP5}, on a $m'\leq
r < 4^n$; par \eqref{m'5}, on a $n_5(m')=n_5(r+4)-1=n_5(r)\geqslant  1$ et,
par \eqref{San}, on voit  que
$\cs(a_n-1)\cap \cs(m')\neq \varnothing$. 
Le lemme \ref{lemPQ} (ii) donne $w \clet w'=M$.

Ainsi, dans tous les cas, les exposants de $U_nP_{r+4}$ sont
domin\'es par $M$ qui est ainsi, par \eqref{P54n+r}, 
l'exposant dominant de $P_k$.
\item
{\bf $k$ et $r$  impairs et $n_5(r)=0$.} On a $n_5(k)=2^{n-1} > 0$.
Nous allons
montrer que le mon\^ome dominant de $P_{k}$ est celui de $U_nP_{r+4}$.
On a $4\notin \cs(r)$, $n_3(r+4)=n_3(r), n_5(r+4)=n_5(r)+1=1$; 
par le lemme \ref{lemP50} (i), l'exposant dominant de $U_n$
est ${a_n-1}$ et, par l'hypoth\`ese de r\'ecurrence, 
on a $P_{r+4}\neq 0$ et l'exposant 
 dominant ${m'}$ de $P_{r+4}$ a
pour code $ [n_3(r),0]$. Par \eqref{degP5}, on a 
$m'\leqslant r < 4^n$, ce qui, par \eqref{San}, implique  $\cs(m')\cap
\cs(a_n-1)=\varnothing$.
Le lemme \ref{lemPQ} (iii) montre que l'exposant dominant de $U_nP_{r+4}$ est  
$$M'=a_n-1+m' \simeq [n_3(r), 2^{n-1}-1]=[n_3(k),n_5(k)-1]$$
et $h(M')=h(k)-1$.

La majoration \eqref{propP523} est encore valable; elle montre que
$h(V_nP_{r+1}) < h(M')$ et qu'ainsi 
les exposants du produit $V_nP_{r+1}$ sont domin\'es par ${M'}$.

Les exposants de $W_nP_r$ sont  aussi domin\'es par ${M'}$~:
comme $n_5(r)=0$, 
l'hypoth\`ese de r\'ecurrence  implique $h(P_r)\leqslant h(r)-3$ et,
par \eqref{hPQ} et \eqref{hWn}, il vient
$$h(W_nP_r)\leqslant h(W_n)+h(P_r)\leqslant 2^{n-1}+h(r)-3=h(k)-3 < h(M').$$
On conclut, par \eqref{P54n+r} et \eqref{hP+Q}, 
que l'exposant dominant de $P_k$ est ${M'}$, ce 
qui, avec le sous-cas (a),  prouve le point (i). 

\item
{\bf $k$ et $r$ pairs.} Nous allons prouver (iii) \`a  
partir de la relation \eqref{P54n+r}.
En utilisant \eqref{hPQ}, \eqref{hWn} et l'hypoth\`ese de r\'ecurrence, 
il vient
\begin{equation}\label{propP524}
h(W_n P_r)\leqslant h(W_n)+h(P_r)\leqslant 2^{n-1}+h(r)-2=h(k)-2.
\end{equation}
En notant que $V_n$ et $P_{r+1}$ sont impairs, on a la majoration
\begin{eqnarray}\label{propP525}
h(V_n P_{r+1})&\leq& h(V_n)+h(P_{r+1})+1 \leqslant 2^{n-1}+h(r+1)-1\notag\\
& &=2^{n-1}+h(r)-1 = h(k)-1.
\end{eqnarray}
De m\^eme, par \eqref{hPQ}, \eqref{hUn}, l'hypoth\`ese de r\'ecurrence
et \eqref{hk+4}, on a
\begin{eqnarray}\label{propP526}
h(U_n P_{r+4})&\leq& h(U_n)+h(P_{r+4}) \leqslant  2^{n-1}-1+h(r+4)-2\notag\\
&\leq& 2^{n-1}+h(r)-2=h(k)-2.
\end{eqnarray}
Finalement, par \eqref{P54n+r} et \eqref{hP+Q}, 
de \eqref{propP524}, \eqref{propP525} et \eqref{propP526},
on d\'eduit 
$$h(P_{k})\leqslant h(k)-1.$$
Mais, par \eqref{hP5mod2}, $h(P_{k})-h(k)$ est congru  \`a  $k$  modulo $2$,
donc est pair, ce qui entra\^ine $h(P_{k})\leqslant h(k)-2$, c'est-\`a-dire
\eqref{hP5k2}.
\end{enumerate}
\item
{\bf Deuxi\`eme cas~: $2\cdot 4^n \leqslant k < 4^{n+1}$.}
On pose $s=k-2\cdot 4^n$ et, par \eqref{P524n+k}, on a
\begin{equation}\label{P524n+s}
P_{k} = P_{2\cdot 4^n+s} = Y_nP_s+x^3U_n^2P_{s+1} +V_n^2P_{s+2} +xU_n^2P_{s+3}.
\end{equation}
Puisque $s\leqslant k-32$, on peut appliquer l'hypoth\`ese de r\'ecurrence \`a
$s, s+1, s+2$ et  $s+3$. On a aussi $s = k-2\cdot 4^n < 2\cdot 4^n$,
$n_3(k)=2^n+n_3(s)$, $n_5(k)=n_5(s)$ et $h(k)=2^n+h(s)$. 

Distinguons trois sous-cas correspondants aux points (i), (ii) et (iii).
\begin{enumerate}
\item
{\bf $k$ et $s$ impairs et $n_5(k)=n_5(s)\geqslant  1$.} Nous
allons montrer que le mon\^ome dominant de $Y_nP_s$ est aussi celui de
$P_k$.

Par le lemme \ref{lemP50} (iv), l'exposant dominant 
de $Y_n$ est ${2^{2n+1}}$. Par
l'hypoyh\`ese de r\'ecurrence, l'exposant 
dominant $m$ de $P_s$ a pour code $[n_3(s),n_5(s)-1]$. Par \eqref{degP5},
on a $m < s < 2^{2n+1}$, ce qui implique $\cs(m)\cap \cs(2^{2n+1})=\varnothing$
et, par le lemme \ref{lemPQ}
(iii), l'exposant dominant  de $Y_nP_s$ est 
\begin{equation}\label{propP527}
M=2\cdot 4^n+m\simeq [2^n+n_3(s),n_5(s)-1]=[n_3(k),n_5(k)-1].
\end{equation}

Nous allons montrer maintenant que les exposants 
de $x^3U_n^2P_{s+1}$, $V_n^2P_{s+2}$ et $xU_n^2P_{s+3}$ sont tous 
domin\'es par $M$. Les majorations ci-dessous s'effectuent 
en utilisant \eqref{hPQ} ainsi que le lemme \ref{lemP50},
l'hypoth\`ese de r\'ecurrence 
et les in\'egalit\'es \eqref{hk+1}, \eqref{hk+2} et \eqref{hk+3}.
On obtient successivement
 \begin{eqnarray}\label{propP528}
h(x^3U_n^2 P_{s+1})&\leq& h(x^3U_n^2)+h(P_{s+1})\leqslant 2^{n}-1+h(s+1)-2\notag\\
&\leq& 2^{n}+h(s)-2 = h(k)-2,
\end{eqnarray}
\begin{eqnarray}\label{propP529}
h(V_n^2 P_{s+2})&\leq& h(V_n^2)+h(P_{s+2}) \leqslant 2^{n}-2+h(s+2) -1\notag\\
&\leq& 2^{n}-2+h(s) = h(k)-2
\end{eqnarray}
et
\begin{eqnarray}\label{propP5210}
h(xU_n^2 P_{s+3})&\leq& h(xU_n^2)+h(P_{s+3}) \leqslant 2^{n}-2+h(s+3)-2\notag\\
&\leq& 2^{n}+h(s)-3 = h(k)-3.
\end{eqnarray}
Comme, par \eqref{propP527}, on a $h(M)= h(k)-1$,
les formules \eqref{P524n+s}, \eqref{propP528}, 
\eqref{propP529} et \eqref{propP5210} montrent que l'exposant dominant 
de $P_{k}$ est $M$ ce qui prouve (i).
\item
{\bf $k$ et $s$  impairs et $n_5(k)=n_5(s)= 0$.}
L'hypoth\`ese de r\'ecurrence donne 
$h(P_s)\leqslant h(s)-3$, impliquant
$$h(Y_nP_s)\leqslant h(Y_n)+h(P_s)=2^n+h(P_s)\leqslant 2^n+h(s)-3=h(k)-3.$$
Par ailleurs, les in\'egalit\'es \eqref{propP528}, \eqref{propP529} et  
\eqref{propP5210} restent valables et, avec \eqref{P524n+s} 
et \eqref{hP+Q}, entra\^inent
$$h(P_{k})\leqslant h(k)-2.$$
Mais, par \eqref{hP5mod2}, $h(P_{k})-h(k)$ est congru modulo $2$ \`a  $k$, 
donc est impair; on obtient ainsi 
$h(P_{k})\leqslant h(k)-3$, ce qui prouve (ii).
\item
{\bf $k$ et $s$ pairs.}
Nous allons prouver (iii) \`a  partir des relations \eqref{P524n+s},
\eqref{hPQ}, du lemme \ref{lemP50}, de l'hypoth\`ese de
r\'ecurrence et des formules \eqref{hk+1}--\eqref{hk+3} qui,
comme $s$ est pair, donnent
$h(s+1)=h(s)$, $h(s+2)\leqslant h(s)+1$ et $h(s+3)\leqslant h(s)+1$.
On obtient successivement
$$h(Y_nP_s)\leqslant h(Y_n)+h(P_s)=2^n+h(P_s)\leqslant 2^n+h(s)-2=h(k)-2,$$
\begin{eqnarray*}
h(x^3U_n^2 P_{s+1})&\leq& h(x^3U_n^2)+h(P_{s+1})+1= 2^{n}+h(P_{s+1})\\
&\leq& 2^{n}+h(s+1)-1= 2^{n}+h(s)-1=h(k)-1,
\end{eqnarray*}
\begin{eqnarray*}
h(V_n^2 P_{s+2})&\leq& h(V_n^2)+h(P_{s+2})= 2^{n}-2+h(P_{s+2})\\
&\leq& 2^{n}-4+h(s+2)\leqslant 2^{n}-3+h(s)=h(k)-3
\end{eqnarray*}
et
\begin{eqnarray*}
h(xU_n^2 P_{s+3})&\leq& h(xU_n^2)+h(P_{s+3})+1= 2^{n}-1+h(P_{s+3})\\
&\leq& 2^{n}-2+h(s+3)\leqslant 2^{n}+h(s)-1=h(k)-1.
\end{eqnarray*}
Par les formules \eqref{P524n+s} et \eqref{hP+Q}, il vient donc
$$h(P_{k})\leqslant h(k)-1$$
mais, par l'argument de parit\'e \eqref{hP5mod2}, $h(P_{k}) - h(k)$
est pair; cela implique $h(P_{k})\leqslant h(k)-2$, ce qui prouve
(iii).\hfill $\Box$
\end{enumerate}
\end{enumerate}

%%%%%%%%%%%%%%%%%%%%%%%%%%%%%%%%%%%%%%%%%%%%%%%%%%%%%%%%%%%%%%%%
%%%%%%%%%%%%%%%%%%%%%%%%%%%%%%%%%%%%%%%%%%%%%%%%%%%%%%%%%%%%%%%
\section{D\'etermination de l'ordre de nilpotence}% \ref{th1}}
\label{parordregk}
%%%%%%%%%%%%%%%%%%%%%%%%%%%%%%%%%%%%%%%%%%%%%%%%%%%%%%%%%%%%%
%%%%%%%%%%%%%%%%%%%%%%%%%%%%%%%%%%%%%%%%%%%%%%%%%%%%%%%%%%%%%%%%%

%%%%%%%%%%%%%%%%%%%%%%%%%%%%%%%%%%%%%%%%%%%%%%%%%%%%%%%%%%%%%%%
\subsection{Calcul de l'ordre de nilpotence}% \ref{th1}}
\label{demth1}
%%%%%%%%%%%%%%%%%%%%%%%%%%%%%%%%%%%%%%%%%%%%%%%%%%%%%%%%%%%%%

Rappelons que $\cf$ est le $\F_2$-espace vectoriel 
de $\F_2[\D]$ engendr\'e par les puissances d'exposant impair de $\D$
(cf. \S~\!\ref{parFi}). 
En utilisant la relation de domination \eqref{ord}, nous \'ecrirons une
forme modulaire $f\in \cf$, $f\neq 0$ sous la forme
\begin{equation}\label{fm1m2}
f=\D^{m_1}+\D^{m_2} \ldots +\D^{m_r} \;\;  \text{ avec }\;\;
m_1 \cget m_2 \cget \ldots \cget m_r.
\end{equation}
On dit que $m_1$ est l'exposant dominant de $f$
et l'on d\'efinit $h(f)$ par  \eqref{hP}, c'est-\`a-dire
\begin{equation}\label{hf}
h(f)=h(m_1)=\max_{1 \leqslant i \leqslant r} h(m_i).
\end{equation}

\begin{theorem}\label{thmg=h}
(i) Pour $p=3$ ou $5$, on a
$$h(T_p|f) \leqslant h(f)-1.$$

(ii) Si $n_3(m_1)\geqslant  1$, on a $T_3|f\neq 0$ et l'exposant dominant de
$T_3|f$ a pour code $[n_3(m_1)-1,n_5(m_1)]$.

(iii) Si $n_5(m_1)\geqslant  1$, on a $T_5|f\neq 0$ et l'exposant dominant de
$T_5|f$ a pour code $[n_3(m_1),n_5(m_1)-1]$.

(iv) On a
\begin{equation}\label{T3T5fD}
T_3^{n_3(m_1)} T_5^{n_5(m_1)} |f=\D.
\end{equation}

(v) La valeur de l'ordre de nilpotence $g(f)$ 
(cf. \S~\!\ref{parindnil}) est donn\'ee par
\begin{equation}\label{gf=hf}
g(f) = h(f)+1.
\end{equation}
\end{theorem}

\begin{dem}
(i) En utilisant la notation $P_k^{(p)}=T_p|\D^k$ (cf. \eqref{Trk3} et
\eqref{Trk5}), on a
$$T_p|f=T_p\left| \left(\sum_{i=1}^r \D^{m_i}\right)\right. =
\sum_{i=1}^r T_p|\D^{m_i}=\sum_{i=1}^r P_{m_i}^{(p)}(\D)$$
ce qui, par \eqref{hP+Q}, entra\^ine
$$h(T_p|f)\leqslant \max_{1\leqslant i\leqslant r} h(P_{m_i}^{(p)}).$$
Par \eqref{hP3k}  (si $p=3$) et 
par \eqref{hP5k}  (si $p=5$) on a
$$h(P_{m_i}^{(p)})\leqslant h(m_i)-1\quad (1\leqslant i \leqslant r)$$
et (i) en d\'ecoule, puisque, par \eqref{hf}, $h(m_i) \leqslant h(f)$.

\mk

(ii) D\'esignons par $\nu_i$ (lorsque $T_3|\D^{m_i}\neq 0$) l'exposant 
dominant de $T_3|\D^{m_i}=P_{m_i}^{(3)}$. Par 
la proposition \ref{propP32ii} (i), on a
$\nu_1\simeq [n_3(m_1)-1,n_5(m_1)]$, et il suffit de d\'emontrer que, pour
$i\geqslant  2$, 
\begin{equation}\label{ming1}
\nu_i \;\clet \;\nu_1.
\end{equation}
Si $h(m_i) < h(m_1)$, par \eqref{hP3k}, il vient 
$$h(\nu_i) \leqslant h(m_i) -1 <  h(m_1) -1=h(\nu_1),$$ 
ce qui prouve \eqref{ming1}.

\noindent
Si $h(m_i) = h(m_1)$, comme $m_i$ est domin\'e par $m_1$, on
a $n_5(m_i) <  n_5(m_1)$ (si $n_5(m_1)=0$, ce cas ne se pr\'esente
pas), $n_3(m_i) > n_3(m_1) \geqslant  1$ et, par la proposition
\ref{propP32ii} (i), $\nu_i \simeq [n_3(m_i)-1,n_5(m_i)]$; 
\eqref{ming1} en r\'esulte.

\mk

(iii) La preuve est identique \`a  celle de (ii), en appliquant
\eqref{hP5k} et la proposition \ref{propP52} (i) au lieu de 
respectivement \eqref{hP3k} et la proposition \ref{propP32ii} (i).

\mk

(iv) Posons $\f=T_3^{n_3(m_1)} T_5^{n_5(m_1)}|f$. 
En appliquant $n_3(m_1)$ fois (ii) et $n_5(m_1)$ fois (iii), on voit
que $\f\neq 0$ et que son exposant dominant $m$ a pour
code $[0,0]$; comme $m$ est impair, on a $m=1$ d'o\`u $\f=\D$, ce qui 
d\'emontre \eqref{T3T5fD}. Notons que \eqref{T3T5fD} implique
\begin{equation}\label{T3T5gf}
g(f) \geqslant  n_3(m_1)+n_5(m_1)+1=h(m_1)+1=h(f)+1.
\end{equation}

(v) Soit $d=\max(m_1,m_2,\ldots,m_r)$ le degr\'e de $f$; on va d\'emontrer 
\eqref{gf=hf} par r\'ecurrence sur le nombre impair $d$.

Si $d=1,3$ ou $5$, \eqref{gf=hf} r\'esulte de \eqref{g1g3}, \eqref{g5} et
de la table \eqref{tablen3+n5}.

Soit $d \geqslant  7$ et supposons \eqref{gf=hf} vraie pour toute 
forme de degr\'e $\leqslant d-2$. 
Pour $d \geqslant  7$, on a $h(d) \geqslant  2$ et la d\'efinition de
l'exposant dominant entra\^ine $h(f)=h(m_1) \geqslant  h(d) \geqslant  2$. Par \eqref{T3T5gf}, on a
$g(f)\geqslant  h(f)+1 \geqslant  3$; donc il existe des nombres premiers impairs 
$p_1,p_2,\ldots, p_s$ avec $s=g(f)-1 \geqslant  2$ et 
\begin{equation}\label{ming2}
T_{p_1} T_{p_2}\ldots T_{p_s}|f \neq 0.
\end{equation}
Posons $\f=T_{p_s}|f$ et calculons $g(\f)$. De \eqref{ming2}, on d\'eduit
$$T_{p_1} T_{p_2}\ldots T_{p_{s-1}}|\f =T_{p_1} T_{p_2}\ldots T_{p_s}|f \neq 0,$$ 
ce qui implique $g(\f) \geqslant  s$. Mais \eqref{gTpf} entra\^ine
$g(\f)=g(T_p|f)\leqslant g(f)-1=s$. On en d\'eduit
\begin{equation}\label{ming4}
g(\f)= s =g(f)-1\geqslant  2.
\end{equation} 
Observons que \eqref{ming2} et $s\geqslant  2$ entra\^inent $\f\neq 0$.
Par \eqref{Tpfk}, le degr\'e de $\f$ est $\leqslant d-2$;
on peut donc appliquer \`a  $\f$ l'hypoth\`ese de r\'ecurrence, ce 
qui donne $g(\f)=h(\f)+1$. En d\'esignant par $j$ l'exposant
dominant de $\f$, avec \eqref{ming4}, il vient
\begin{equation}\label{ming5}
g(\f)=h(\f)+1=h(j)+1=s\geqslant  2.
\end{equation}
Soit $[u,v]$ le code de $j$, avec $u\geqslant  0$, $v\geqslant  0$ et $u+v=s-1$.
En appliquant (iv) \`a $\f$ et en posant 
$q_1=q_2=\ldots =q_u=3$ et
$q_{u+1}=q_{u+2}=\ldots = q_{u+v}=5$, il vient
$$T_{q_1} T_{q_2}\ldots T_{q_{s-1}}|\f = 
T_{q_1} T_{q_2}\ldots T_{q_{s-1}}T_{p_s}|f = \D.$$
Posons $\psi=T_{q_{s-1}}|f$; on a
$$T_{q_1} T_{q_2}\ldots T_{q_{s-2}}T_{p_s}|\psi = 
T_{q_1} T_{q_2}\ldots T_{q_{s-1}}T_{p_s}|f = \D.$$ 
Cette formule montre que $g(\psi)\geqslant  s$. 
Mais \eqref{gTpf} entra\^ine
$g(\psi)=g(T_{q_{s-1}}|f)\leqslant g(f)-1=s$ et $g(\psi)=s$.

Par \eqref{Tpfk}, le degr\'e de $\psi$ est  $\leqslant d-2$;
et l'hypoth\`ese de r\'ecurrence 
donne $g(\psi)=h(\psi)+1$. On a ainsi
\begin{equation}\label{ming6}
g(\psi)= s =g(f)-1=h(\psi)+1.
\end{equation}
Par (i), on a $h(T_{q_{s-1}}|f) \leqslant h(f)-1$, d'o\`u, par \eqref{ming6}  
$$s-1=g(f)-2 = h(\psi)=h(T_{q_{s-1}}|f) \leqslant h(f)-1$$
ce qui implique $g(f) \leqslant h(f)+1$; vu \eqref{T3T5gf},
cela entra\^ine \eqref{gf=hf}.
\end{dem}

\begin{coro}\label{coroT3T5D}
Soit $f\in\cf$, $f\neq 0$. Si $T_3|f=T_5|f=0$, alors $f=\D$.
\end{coro}

\begin{dem}
En effet, d'apr\`es (iv), on a $n_3(m_1)=n_5(m_1)=0$, 
d'o\`u $m_1=1$ et $f=\D$.
\end{dem}

\begin{coro}\label{corothm2}
Soit $f\in\cf$, $f\neq 0$, et $p$ un nombre premier v\'erifiant 
$p\equiv \pm 1 \pmod{8}$. Alors, on a 
\begin{equation}\label{corothm21}
g(T_p|f)\leqslant g(f)-2.
\end{equation}
\end{coro}

\begin{dem}
Consid\'erons d'abord le cas o\`u $f=\D^k$ 
(avec $k$ impair, $k\geqslant  1$). 
\eqref{gf=hf} implique $g(\D^k)=h(\D^k)+1=h(k)+1\geqslant  1$. 

Si $T_p|\D^k=0$, par \eqref{g0} on a $g(T_p|\D^k)=-\iy$ et
\eqref{corothm21} est v\'erifi\'ee.

Si $T_p|\D^k \neq 0$, on \'ecrit 
$T_p|\D^k=\D^{m_1}+\D^{m_2}+ \ldots +\D^{m_r}$ et, par \eqref{Tpfk}, 
pour $1\leqslant i \leqslant r$, on a
$m_i\equiv pk\equiv \pm k \pmod{8}$.
Par la table \eqref{tablen3n5h}, il s'ensuit que 
$h(m_i)\equiv h(k)\pmod{2}$. Ainsi,
par \eqref{hf}, on a
\begin{equation}\label{corothm22}
h(T_p|\D^k)\equiv h(k) \pmod{2}.
\end{equation}
Ensuite, en appliquant \eqref{gf=hf} \`a $T_p|\D^k$ et \`a $\D^k$, on
obtient
\begin{equation}\label{corothm23}
g(T_p|\D^k)=h(T_p|\D^k)+1 \quad \text{ et } \quad g(\D^k)=h(\D^k)+1=h(k)+1.
\end{equation}
On d\'eduit de \eqref{corothm22} et \eqref{corothm23}
\begin{equation}\label{corothm24}
g(T_p|\D^k)\equiv g(\D^k) \pmod{2}.
\end{equation}
Mais, par \eqref{gTpf}, on a $g(T_p|\D^k)\leqslant g(\D^k)-1$, ce qui avec 
\eqref{corothm24}, prouve 
\begin{equation}\label{corothm241}
g(T_p|\D^k)\leqslant g(\D^k)-2
\end{equation}
c'est-\`a-dire \eqref{corothm21} lorsque $f$ est un mon\^ome en $\D$.
Notons que \eqref{corothm23} et \eqref{corothm241} entra\^inent
\begin{equation}\label{corothm25}
h(T_p|\D^k)\leqslant h(\D^k)-2=h(k)-2.
\end{equation}

Soit maintenant $f=\D^{k_1}+\D^{k_2}=\ldots + \D^{k_s}$ avec $s\geqslant  2$ et
$k_1 > k_2  > \ldots > k_s$. Il vient $T_p|f=T_p|\D^{k_1}+T_p|\D^{k_2}+\ldots +
T_p|\D^{k_s}$ et de \eqref{hP+Q} et \eqref{corothm25}, il r\'esulte
\begin{eqnarray*}
h(T_p|f) &\leq& \max(h(T_p|\D^{k_1}),h(T_p|\D^{k_2}),\ldots ,h(T_p|\D^{k_s}))\\
&\leq& \max(h(k_1)-2, h(k_2)-2,\ldots, h(k_s) -2)=h(f)-2,
\end{eqnarray*}
ce qui, avec les relations $g(T_p|f)=h(T_p|f)+1$ et $g(f)=h(f)+1$ 
fournies par \eqref{gf=hf}, termine la d\'emonstration 
de \eqref{corothm21}.
\end{dem}

%%%%%%%%%%%%%%%%%%%%%%%%%%%%%%%%%%%%%%%%%%%%
\subsection{Une autre preuve de la nilpotence de $T_p$}
\label{parpreuvenil}
%%%%%%%%%%%%%%%%%%%%%%%%%%%%%%%%%%%%%%%%%%%%%

Soit $p$ un nombre premier congru \`a $3,5$ ou $7$ modulo 2. Nous avons
vu au \S\ref{parnil} que la nilpotence de $T_p$ r\'esultait de
\eqref{Tpfk1} et de \eqref{muj0}.  

Soit maintenant $p\equiv 1 \bmod 8$. 
Par \eqref{Tpfk1}, $T_p|\D$ est \'egal \`a $0$ ou \`a $\D$. 
Mais $T_p(\D)=\D$ entrainerait que
dans le d\'eveloppement en $q$ de $\D$ le coefficient de $q^p$ soit
\'egal \`a $1$,  ce qui, par \eqref{Delta}, n'est pas
vrai. On a donc  
\begin{equation}\label{TpDk=0}
T_p|\D= 0.
\end{equation}
Supposons qu'il existe $k$ impair tel que $T_p|\D^k$
soit un polyn\^ome de degr\'e $k$ en $\D$
et soit $k_0$ le plus petit tel $k$. 
Notons que \eqref{TpDk=0} implique $k_0 \ge 3$. On peut \'ecrire
\begin{equation}\label{TpDk_0}
T_p|\D^{k_0} =\D^{k_0}+ \sum_{j \leqslant k_0-2} \e_{k_0,j} \;\D^j,
\quad \e_{k_0,j}\in \F_2
\end{equation}
et, pour $k$ impair $< k_0$,
\begin{equation}\label{TpDkk}
T_p|\D^{k} =\sum_{j \leqslant  k-2} \e_{k,j} \;\D^j,
\quad \e_{k,j}\in \F_2.
\end{equation}
Soit $\cf^{(k_0)}$ le $\F_2$-espace vectoriel de 
base $\D,\D^3,\ldots,\D^{(k_0)}$; 
la matrice de la restriction de l'op\'erateur $T_p$ \`a  
$\cf^{(k_0)}$ dans cette base est triangulaire, 
ses valeurs propres sont $0$ et $1$ et il existe une forme
$\f\in \cf^{(k_0)}$ de degr\'e $k_0$ v\'erifiant 
\begin{equation}\label{Tpphi=phi}
T_p|\f= \f.
\end{equation}
Remarquons que $\f\neq \D$ puisque le
degr\'e $k_0$ de $\f$ est $\geqslant  3$. Par le corollaire \ref{coroT3T5D}, 
il existe $\ell\in \{3,5\}$ tel que
\begin{equation}\label{Tlphi<>0}
T_\ell|\f\neq 0.
\end{equation}
Ensuite, on fixe $r > k_0$ et l'on consid\`ere 
l'action de l'op\'erateur
$$T_\ell (T_p)^r=(T_p)^r T_\ell$$
sur $\D^{k_0}$. Par \eqref{Tpphi=phi} et \eqref{Tlphi<>0},
il vient
\begin{equation}\label{TlTpr<>0}
T_\ell (T_p)^r|\f=T_\ell|\f \neq 0.
\end{equation}
D'autre part, puisque $T_3$ et $T_5$ sont nilpotents, on a 
$$\text{degr\'e } T_\ell|\f \leqslant k_0-2$$
ce qui implique, par \eqref{TpDkk}
\begin{equation*}
(T_p)^r \,T_\ell|\f =0,
\end{equation*}
en contradiction avec \eqref{TlTpr<>0}, ce qui 
ach\`eve la preuve de la nilpotence de l'op\'erateur de Hecke
$T_p$ modulo $2$.

%%%%%%%%%%%%%%%%%%%%%%%%%%%%%%%%%%%%%%%%%%%%
\subsection{Estimation de $g(k)=g(\D^k)$}
\label{parestimgk}
%%%%%%%%%%%%%%%%%%%%%%%%%%%%%%%%%%%%%%%%%%%%%
 
\begin{prop}\label{propg}
Soit $k$ impair $\geqslant  1$. On a
\begin{equation}\label{ineggk}
\frac{\sqrt k}{2} < 1+\frac{\sqrt{k-1}}{2} \leqslant g(k) \leqslant 
\frac 32 \sqrt{k+1}-1 < \frac 32 \sqrt k.
\end{equation}
On a $g(k)=1+\frac{\sqrt{k-1}}{2}$ si et seulement si $k$ est de la forme
$$k=1 \quad \text{ ou } \quad k=4^a+1\simeq [0,2^{a-1}]\quad 
\text{ avec } a\geqslant  1.$$
On a $g(k)=\frac 32 \sqrt{k+1}-1$ si et seulement si $k$ est de la forme
\begin{equation}\label{k1}
k=1+2+4+8+\ldots + 2^{2a-1}=4^{a}-1\simeq [2^a-1,2^{a-1}-1].
\end{equation}
\end{prop}
 
\begin{dem}
Soit $r$ le plus petit nombre tel que $k\leqslant 4^r-1$.
On pose
$$k=1+\sum_{i=1}^{2r-1} \b_i 2^i,\quad \text{ avec } \b_i\in \{0,1\} \quad
\text{ et } \quad \max(\b_{2r-1},\b_{2r-2})=1.$$
Soit $i\geqslant  1$; on pose (cf \eqref{h2i})
$$U_i=h(2^i)=
\begin{cases}
n_5(2^i)=2^{\frac{i-2}{2}} & \text{ si $i$ est pair}\\
n_3(2^i)=2^{\frac{i-1}{2}} & \text{ si $i$ est impair},
\end{cases}
$$
$V_0=0$ et
$$V_i=\sum_{1\leqslant j \leqslant i} U_j=
\begin{cases}
2^{\frac{i+2}{2}} -2& \text{ si $i$ est pair}\\
3\cdot2^{\frac{i-1}{2}} -2& \text{ si $i$ est impair}.
\end{cases}
$$
Par \eqref{h}, on a ainsi
$$h(k)=\sum_{i=1}^{2r-1} \b_i U_i \leqslant V_{2r-1},$$
et, par le th\'eor\`eme \ref{thmg=h},
$$g(k)=g(\D^k)=h(k)+1.$$
{\bf Minoration de $g(k)$.} En observant que $\b_i^2=\b_i$, il vient 
\begin{eqnarray}
(g(k)-1)^2=h(k)^2=\left(\sum_{i=1}^{2r-1} \b_i U_i \right)^2 &\geq&
\sum_{i=1}^{2r-1} \b_i^2 U_i^2\label{mingk1}\\ 
&\geq& \sum_{i=1}^{2r-1} \b_i 2^{i-2}=\frac{k-1}{4}\cdot\label{mingk2}
\end{eqnarray}
L'in\'egalit\'e large \eqref{mingk1} est une \'egalit\'e si et seulement si
au plus un chiffre binaire $\b_i$ de $k-1$ est non nul. Il n'y a 
\'egalit\'e dans \eqref{mingk2} que si ce chifre $\b_i$ a  un indice $i$ pair. 

\ni
{\bf Majoration de $g(k)$.} On a 
\begin{eqnarray}\label{majgk}
(1+g(k))^2 &=& (2+h(k))^2=\left(2+\sum_{i=1}^{2r-1} \b_i U_i \right)^2 \notag\\
&=& 4+ \sum_{i=1}^{2r-1}\b_iU_i\left(4+U_i+
2\sum_{1\leqslant j \leqslant i-1}\b_j U_j\right)\notag\\
&\leq& S \;\stackrel{def}{=\!=} \;4+ \sum_{i=1}^{2r-1}\b_iU_i\left(4+U_i+
2V_{i-1}\right)\notag
\end{eqnarray}
en utilisant l'in\'egalit\'e $\b_i\sum_{1\leqslant j \leqslant i-1}\b_j U_j\leqslant \b_iV_{i-1}$.
Notons que cette in\'egalit\'e n'est une \'egalit\'e que si $\b_i=0$ ou 
$\b_1=\b_2=\ldots=\b_i=1$. 
Cela entra\^ine, puisque l'un des deux 
chiffres $\b_{2r-1}$ ou $\b_{2r-2}$ est non nul
\begin{equation}\label{S=1+gk}
(1+g(k))^2= S \quad \Llr \quad \b_1=\b_2=\ldots=\b_{2r-2}=1.
\end{equation}
Calculons maintenant $S$, en s\'eparant les indices pairs et impairs.
\begin{eqnarray*}
S&=&4+\b_1U_1(4+U_1)+ \sum_{\ell=1}^{r-1}\Bigl(\b_{2\ell} U_{2\ell} \left(4+U_{2\ell}+
2 V_{2\ell-1}\right)\\
& & \hspace{20mm} + \b_{2\ell +1} U_{2\ell+1} \left(4+U_{2\ell+1}+
2 V_{2\ell-1}+2\b_{2\ell} U_{2\ell}\right)\Bigr)\notag\\
&=&  4+5\b_1+\sum_{\ell=1}^{r-1} 
2^{2\ell}\left(\frac 74\b_{2\ell}+4\b_{2\ell+1}+\b_{2\ell}\b_{2\ell+1}\right)\notag.
\end{eqnarray*}
Il s'ensuit que
\begin{eqnarray}\label{majgk1}
 \frac 94(k+1)-S &=& \frac 94(2+2\b_1)-(4+5\b_1)\notag\\
+ \sum_{\ell=1}^{r-1} 2^{2\ell} \hspace{-7mm}& & \hspace{-2mm}
\left(\frac 94(\b_{2\ell}+2\b_{2\ell+1})-
\frac 74\b_{2\ell}-4\b_{2\ell+1}-\b_{2\ell}\b_{2\ell+1}\right)\notag\\
&=& \frac 12(1-\b_1)+\sum_{\ell=1}^{r-1} 
2^{2\ell-1}(\b_{2\ell}-\b_{2\ell+1})^2\geqslant  0.
\end{eqnarray}
Par \eqref{majgk} et \eqref{majgk1}, pour que l'on ait
$(1+g(k))^2=\frac 94 (k+1)$ il faut et il suffit que l'on ait
$\b_1=\b_2=\ldots=\b_{2r-1}=1$, c'est-\`a-dire $k=4^r-1$.
\end{dem}

\begin{coro}\label{coroordnilgf}
Soit $f\in \F_2[\D]$ une forme modulaire modulo $2$ de degr\'e $d\geqslant  1$
en $\D$. On a
$$g(f) < \frac 32 \sqrt{d}.$$ 
\end{coro}

\begin{dem}
Consid\'erons d'abord le cas $f\in \cf$, que l'on \'ecrit sous la
forme \eqref{fm1m2}.  Par \eqref{gf1f2} et \eqref{ineggk}, on a
$$g(f) \leqslant \max (g(m_1), \ldots, g(m_r)) <
\frac 32 \max(\sqrt{m_1},\ldots,\sqrt{m_r})
=\frac 32 \sqrt{d}.$$
Supposons maintenant $f$ parabolique. Par \eqref{fs},
$f=\sum_{s=0}^S f_s^{2^s}$ avec $f_s\in \cf$. Mais, pour $p$ impair,
on a $T_p| f_s^{2^s}=\left(T_p|f_s\right)^{2^s}$, d'où l'on d\'eduit $g(
f_s^{2^s})=g(f_s)$ et, par \eqref{gf1f2}, en notant $d_s$ le degr\'e de
$f_s$ (si $f_s=0, d_s=-\iy$), on a
\begin{eqnarray*}
g(f) \hspace{-2mm}&\leq&
\max\left((g\left(f_0^{2^0}\right),g\left(f_1^{2^1}\right),\ldots,
g\left(f_S^{2^S}\right)\right)=\max(g(f_0),g(f_1),\ldots,g(f_S))\\
&<& \frac 32 \max(d_0,d_1,\ldots, d_S) \leqslant \frac 32 d.
\end{eqnarray*}
Enfin, si $f$ est non parabolique, par \eqref{1+fs}, $\f=f-1$ est
parabolique de degr\'e $d\geqslant  1$. Par \eqref{ga2}, $T_p|1=0$ pour tout
$p$ premier impair et $g(1) =1$. On a donc, par \eqref{gf1f2},
 $\ds g(f) =g(1+\f)\leqslant \max(g(1),g(\f)) <\frac 32 \sqrt d$.
\end{dem}

\begin{prop}\label{propn3}
Soit $k$ impair $\geqslant  1$. On a
\begin{equation}\label{inegn3}
0\leqslant n_3(k) \leqslant \sqrt{\frac{3k-1}{2}}-1 <\sqrt{\frac{3k}{2}}\cdot
\end{equation}
On a $n_3(k)=\sqrt{\frac{3k-1}{2}}-1$ si et seulement si $k$ est de la forme
\begin{equation}\label{n3k1}
k=1+2+8+32+\ldots +2\cdot 4^{a-1}=1+2\frac{4^a-1}{3}\simeq [2^a-1,0].
\end{equation}
\end{prop}

\begin{dem}
On a $n_3(k)=\sum_{i=0}^\iy \b_{2i+1} 2^i$ (cf. \S\ref{parn3n5h}) et
\begin{eqnarray*}\label{majn5}
(n_3(k)+1)^2 &=& 1+ 2 \sum_{i=0}^\iy \b_{2i+1} 2^{i}+
\sum_{i=0}^\iy \b_{2i+1}^2 2^{2i}\\
& & \hspace{8mm} +2\sum_{i=0}^\iy
\b_{2i+1} 2^{i} \left(\sum_{0\leqslant j < i} \b_{2j+1} 2^j\right)\\
&\leq& 1+ 2 \sum_{i=0}^\iy \b_{2i+1} 2^{i}+
\sum_{i=0}^\iy \b_{2i+1} 2^{2i}+2\sum_{i=0}^\iy
\b_{2i+1} 2^{i} ( 2^i-1)\\
&=& 3 \sum_{i=0}^\iy \b_{2i+1} 2^{2i} +1=\frac 32 \sum_{i=0}^\iy
\b_{2i+1} 2^{2i+1} +1 \leqslant 1+\frac 32 (k-1).
\end{eqnarray*}
Pour que l'on ait dans le calcul ci-dessus $(n_3(k)+1)^2=1+\frac 32
(k-1)$, il faut et il suffit que l'on ait $\b_i=1$ pour tout $i$ impair et
$\b_i=0$ pour tout $i$ pair $\geqslant  2$.
\end{dem}

\begin{prop}\label{propn5}
Soit $k$ impair $\geqslant  1$. On a
\begin{equation}\label{inegn5}
0\leqslant n_5(k) \leqslant \sqrt{\frac{3k+1}{4}}-1 <\sqrt{\frac{3k}{4}}\cdot
\end{equation}
On a $n_5(k)=\sqrt{\frac{3k+1}{4}}-1$ si et seulement si
 $k$ est de la forme
\begin{equation}\label{n5k1}
k=1+4+16+64+\ldots +4^{a-1}=\frac{4^a-1}{3}\simeq [0,2^{a-1}-1].
\end{equation}
\end{prop}

\begin{dem}
La d\'emonstration est similaire \`a celle de la proposition \ref{propn3}.
\end{dem}

\def\refname{R\'ef\'erences}

\vspace{2cm}

\noindent
Jean-Louis NICOLAS, 

\noindent
Universit\'e de Lyon, CNRS, Universit\'e Lyon 1,

\noindent
Institut Camille Jordan, Math\'ematiques, 

\noindent
43 Bd. du 11 Novembre 1918, 

\noindent
F-69622 Villeurbanne Cedex, France. 

\mk

\noindent
{\tt nicolas@math.univ-lyon1.fr} \hfill{\tt jpserre691@gmail.com}  

\noindent
\tt{http://math.univ-lyon1.fr/$\sim$nicolas/}

\end{document}